\documentclass[10pt,letterpaper]{amsart}
\usepackage{amscd}
\usepackage{amsfonts}
\usepackage{amsmath}
\usepackage{amsrefs}
\usepackage{amssymb}
\usepackage{amsthm}
\usepackage{array}
\usepackage{bm}
\usepackage{dcpic,pictexwd}
\usepackage{fancyhdr}
\usepackage{indentfirst}
\usepackage{ifthen}
\usepackage{accents}
\numberwithin{equation}{section}
\newtheorem{Theorem}{Theorem}[section]

\newtheorem{Lemma}[Theorem]{Lemma}
\newtheorem{Corollary}[Theorem]{Corollary}

\theoremstyle{definition}
\newtheorem{Example}[Theorem]{Example}
\newtheorem{Remark}[Theorem]{Remark}
\newfont{\deffont}{cmbxti10}
\newfont{\german}{eufm10}
\newfont{\mymath}{cmr12}
\newcommand{\CalG}{\mathcal{G}}
\newcommand{\CalSigma}{\mathcal{H}}

\newcommand\lieh{\mathfrak{h}}

\newcommand\hook{\mathbin{\raise0.5pt\hbox{\hbox{{\vbox{\hrule height.4pt width6pt depth0pt}}}\vrule height6pt width.4pt depth0pt}\,}}

\newcommand\cTM{T^*\kern-2ptM}

\newcommand\Cau{\text{\rm A}}

\newcommand\vess{\text{\german v\german e\german s\german s}}

\newcommand\thetaX{\theta_{\kern -1 pt X}}

\newcommand\ann{\text{\rm ann}}
\newcommand\CalPf{\mathcal{P}\text{\it \kern -.3pt f}}
\newcommand\Real{\text{\bf  R}}

\newcommand\TM{T\kern -2pt M}

\newcommand\real{\text{\bf R}}
\newcommand\mycap{\hbox{\ $\rlap{\kern -.3pt $\cap$}\raise.8pt\hbox{$\scriptstyle+$}$\ } }

\newcommand\Ao{{\kern-2.3pt}\stackrel{\scriptscriptstyle o}{A}{}{\kern-2.3pt}}
\newcommand\Bo{{\kern-2.3pt}\stackrel{\scriptscriptstyle o}{B}{}{\kern-2.3pt}}
\newcommand\Ko{{\kern-2.3pt}\stackrel{\scriptscriptstyle o}{K}{}{\kern-2.3pt}}
\newcommand\Co{{\kern-2.3pt}\stackrel{\scriptscriptstyle o}{C}{}{\kern-2.3pt}}
\newcommand\Qo{{\kern-2.3pt}\stackrel{\scriptscriptstyle o}{Q}{}{\kern-2.3pt}}
\newcommand\Mo{{\kern-2.3pt}\stackrel{\scriptscriptstyle o}{M}{}{\kern-2.3pt}}

\newcommand\Xo{{\stackrel{\scriptscriptstyle o}{X}}{\kern-1.3pt}}
\newcommand\Yo{{\stackrel{\scriptscriptstyle o}{Y}}{\kern-1.3pt}}

\DeclareMathOperator{\rank}{rank}

\DeclareMathOperator{\spn}{span}

\newboolean{proofmode}

\newcommand{\StTag}[1]{ \label{st:#1}
\ifthenelse{\boolean{proofmode}}{\ \marginpar{\quad\scriptsize st:#1} }{}      }

\newcommand{\EqTag}[1]{
\ifthenelse{\boolean{proofmode}}
{ {\label{eq:#1}}
  \stepcounter{equation}
  \tag{\theequation \rlap{\kern 23 pt{\scriptsize eq:#1}}}
}
{\label{eq:#1}}
 }

\newcommand{\EqRef}[1]{\eqref{eq:#1}}
\newcommand{\StRef}[1]{\ref{st:#1}}
\newcolumntype{C}{>\scriptstyle>{$}c <{$} }
\newcolumntype{L}{>\scriptstyle >{$} l <{$} }
\newcommand\CalA{\mathcal{A}}
\newcommand\CalB{\mathcal{B}}
\newcommand\CalC{\mathcal{C}}
\newcommand\CalE{\mathcal{E}}
\newcommand\CalI{\mathcal{I}}

\newcommand\CalJ{\mathcal{J}}
\newcommand\CalK{\mathcal{K}}

\newcommand\CalH{\mathcal{H}}

\newcommand\CalS{\mathcal{S}}

\newcommand\CalV{\mathcal{V}}

\newcommand\barI{ \bar I}

\newcommand\sd{\mathbin{ \raise0.0pt\hbox{ \vrule height5pt width.4pt depth0pt}\!\times}}

\newcommand\barC{
	\hbox{\kern 2.3 true pt
	\vbox{\hrule width 6.5  true pt height .3 true pt \kern .9 true pt
	\hbox{\kern -0.8 true pt $C$}}}}

\newcommand\barM{
	\hbox{\kern 2.3 true pt
	\vbox{\hrule width 8.5  true pt height .3 true pt \kern .9 true pt
	\hbox{\kern -2.3 true pt $M$}}}}
\newcommand\barU{
	\hbox{\kern .8 true pt
	\vbox{\hrule width 6.5  true pt height .3 true pt \kern .9 true pt
	\hbox{\kern -.8 true pt $U$}}}}

\newcommand\barCalI{
	\hbox{\kern 4.3 true pt
	\vbox{\hrule width 6.5  true pt height .3 true pt \kern .9 true pt
	\hbox{ \kern -4.3 true pt  $\CalI$}}}}
\newcommand\barXi{
	\hbox{\kern 1 true pt
	\vbox{\hrule width 6.5  true pt height .3 true pt \kern .9 true pt
	\hbox{\kern 1 true pt $\Xi$}}}}

\newcommand\barDelta {\bar \Delta}
\newcommand\mytilde {\sim }

\newcommand\Largehat{\smash{\raise -7.5 pt \hbox{\rm\Large\^{}}}}
\newcommand\LARGEhat{\smash{\raise -7.5 pt \hbox{\rm\LARGE\^{}}}}
\newcommand\hugehat{\smash{\raise -7.5 pt \hbox{\rm\huge\^{}}}}
\newcommand\Hugehat{\smash{\raise -7.5 pt \hbox{\rm\Huge\^{}}}}

\newcommand\Largecheck{\smash{\raise -7.5 pt \hbox{\rm\Large\v{}}}}
\newcommand\LARGEcheck{\smash{\raise -7.5 pt \hbox{\rm\LARGE\v{}}}}
\newcommand\hugecheck{\smash{\raise -7.5 pt \hbox{\rm\huge\v{}}}}
\newcommand\Hugecheck{\smash{\raise -7.5 pt \hbox{\rm\Huge\v{}}}}

\newcommand\hV{{\accentset{\LARGEhat}{V}}}

\newcommand\cV{{\accentset{\LARGEcheck}{V}}}

\newcommand\hCalK{{\accentset{\LARGEhat}{{\mathcal K}}}}

\newcommand\hCalV{{\accentset{\LARGEhat}{{\mathcal V}}}}

\newcommand\cCalV{{\accentset{\LARGEcheck}{{\mathcal V}}}}

\newcommand\bfG{\boldsymbol{G}}

\newcommand\bfq{\mathbf{q}}
\newcommand\bfp{\mathbf{p}}

\newcommand\bfGamma{\boldsymbol{\Gamma}}

\newcommand\bfmu{{\boldsymbol{\mu}}}

\newcommand\bfthetaX{{\boldsymbol{\theta_{\kern -1 pt X}}}}

\newcommand\Rtheta{{ \raise 1pt \hbox{$\scriptstyle {\boldsymbol{\theta}}$}}}
\newcommand\Ltheta{{ \lower 1pt \hbox{$\scriptstyle {\boldsymbol{\theta}}$}}}

\newcommand\Rsigma{{ \raise 1pt \hbox{$\scriptstyle \sigma$}}}

\newcommand\Reta{{ \raise 1pt \hbox{$\scriptstyle \eta$}}}

\newcommand\vecthsigma{\partial_{{\displaystyle \hat {\raise 1.3pt \hbox{$\scriptstyle \sigma$}}}^a}}
\newcommand\vectcsigma{\partial_{{\displaystyle \check {\raise 1.3pt \hbox{$\scriptstyle \sigma$}}}^\alpha}}

\setboolean{proofmode}{false}

\newcommand\alg{\text{alg}}

\def\Obj(#1, #2)[#3]#4{\obj(#1, #2)[#3]{#4}}
\def\beginDC#1[#2]{\begindc{#1}[#2]}

\begin{document}

\title{Symmetry Reduction of Exterior Differential Systems  and B\"acklund Transformations for PDE in the Plane.}

\author{I.M. Anderson}
\address{Department of Mathematics and Statistics, Utah State University, Logan Utah, 84322}
\curraddr{Department of Mathematics and Statistics, Utah State University, Logan Utah, 84322}
\email{Ian.Anderson@usu.edu}
\thanks{Supported in part by NSF Grant \#0713830.}

\author{M.E. Fels}
\address{Department of Mathematics and Statistics, Utah State University, Logan Utah, 84322}
\curraddr{Department of Mathematics and Statistics, Utah State University, Logan Utah, 84322}
\email{Mark.Fels@usu.edu}


\date{\today}



\begin{abstract} We approach the construction of B\"acklund transformations for Darboux integrable hyperbolic partial differential equations in the plane through the reduction of exterior differential systems.
\end{abstract}

\pagestyle{plain}
\maketitle

\section{Introduction}  Let $\CalI$ be an exterior  differential system (EDS) on a manifold $M$ and  let $\bfp: M \to N$ be a smooth submersion. We define the  reduced differential system $ \CalI/\bfp$ on $N$ by
\begin{equation}
	 \CalI/\bfp = \{\, \theta \in \Omega^*(N) \, |\, \bfp^*(\theta ) \in \CalI \, \} .
\EqTag{SQ}
\end{equation}
In the special case where $G$ is a symmetry group of $\CalI$  which acts regularly on $M$, we shall write $\CalI/G$ in place of $\CalI/\bfq_G$,
where $\bfq_G : M \to M/G$ is the canonical projection to the space of orbits $M/G$.  Similarly if $\Gamma_G$ a Lie algebra of infinitesimal symmetries of  $\CalI$ which is regular on $M$,  we shall write $\CalI/\Gamma_G$ in place of $\CalI/\bfq_{\Gamma_G}$, where $\bfq_{\Gamma_G} : M \to M/\Gamma_G$ is the canonical projection to the leaf space $M/\Gamma_G$.

Two  exterior differential systems (EDS) $\CalI_1$ and $\CalI_2$, defined on manifolds $M_1$ and $M_2$ are said to be related by a B\"acklund transformation
if there exists an EDS  $\CalB$ on a manifold $N$ which serves simultaneously as an integrable extensions for both  $\CalI_1$ and $\CalI_2$.
Precisely there are maps  $\pi_1 :N \to M_1$ and  $\pi_2 : N \to M_2$  where
\begin{equation}
\beginDC{\commdiag}[3]
\Obj(0, 14)[B]{$\CalB$}
\Obj(-12, 0)[I1]{$\CalI_1$}
\Obj(12, 0)[I2]{$\CalI_2$}
\mor{B}{I1}{$\pi_1$}[\atright, \solidarrow]
\mor{B}{I2}{$\pi_2$}[\atleft, \solidarrow]
\enddc
\EqTag{BT0}
\end{equation}
and $\CalB$ is the B\"acklund transformation.

B\"acklund transformations  \EqRef{BT0} can then be constructed using EDS reduction  \EqRef{SQ} as follows.
\begin{Theorem}
\StTag{ThIntro1}
 Let $\CalI$ be a Pfaffian system on $M$ and with symmetry groups $G_1$ and $G_2$. Let $H$ be a subgroup of $G_1$ and $G_2$, and assume that
\par\smallskip
\noindent
{\bf [i]}
	$M/H$, $M/G_1$ and $M/G_2$ are smooth manifolds with smooth quotient maps  $\bfq_H : M \to M/H$, $\bfq_{G_1} : M \to M/G_1$, and $\bfq_{G_2} : M \to M/G_2$, then

\begin{equation}
\beginDC{\commdiag}[3]
\Obj(0, 32)[I]{$\CalI$}
\Obj(0, 12)[H]{$\CalI /H$}
\Obj(-20, 0)[I1]{$\CalI /G_1$}
\Obj(20,  0)[I2]{$\CalI /G_2$}
\mor{I}{H}{$\bfq_H$}[\atleft, \solidarrow]
\mor{H}{I1}{$\bfp_1$}[\atleft, \solidarrow]
\mor{H}{I2}{$\bfp_2$}[\atright, \solidarrow]
\mor{I}{I1}{$\bfq_{G_1}$}[\atright, \solidarrow]
\mor{I}{I2}{$\bfq_{G_2.}$}[\atleft, \solidarrow]
\enddc
\EqTag{Intro1}
\end{equation}
is a commutative diagram of EDS, where $\bfp_1$ and $\bfp_2$ are the orbit maps defined by
\begin{equation}
\bfp_1(Hx) = G_1x , \quad {\text and} \quad \bfp_2(Hx)=G_2x \quad {\text for \ all } \ x \in M.
\EqTag{IntExt140}
\end{equation}

\noindent
{\bf [ii]} If the action of $G_1 $ and $G_2$ are transverse to $\CalI$, then $\bfq_{G_1}$, $\bfq_{G_2}$, $ \bfp_1$ and $\bfp_2$ are all integrable extensions, and $\CalI/H$ defines a B\"acklund transformation between
$\CalI /G_1$  and  $\CalI /G_2$.
\end{Theorem}

In this article we will apply Theorem \StRef{ThIntro1} to PDE in the plane. Let  $M \subset J^2(\real^2,\real)$ be 
\begin{equation}
M= \{\ p \in J^2(\real^2, \real) \ | \ F(x,y,u,u_x,u_y, u_{xx},u_{xy},u_{yy})=0 \ \}
\EqTag{DEFM}
\end{equation}
where $F$ is a hyperbolic PDE in the plane and let $\CalI_2$ be the standard rank 3 hyperbolic Pfaffian system determined by the restriction of the contact system on $J^2(\real^2,\real)$ \cite{gardner-kamran:1993a} to the level set $M$ in equation \EqRef{DEFM}. Associated with the hyperbolic system  $\CalI_2$ are a pair of rank $5$ Pfaffian systems $\hV$ and $\cV$ called the characteristic \cite{gardner-kamran:1993a} or singular Pfaffian system \cite{anderson-fels-vassiliou:2009a}.  A function $f\in C^\infty(M)$ is called a Darboux invariant if $df \in\hV$ or $df\in \cV$, while 
the rank of the integrable subsystems $\hV^\infty \subset \hV$ and $\cV^\infty\subset \cV$ determine the number of independent invariant.  The differential system $\CalI_2$ is said to be Darboux integrable (and not Monge integrable)  if $\rank \hV = \rank \cV = 2$. Our first theorem constructs B\"acklund transformations between generic PDE satisfying the Darboux integrability condition and the wave equation given
as an EDS by $\CalI_{s=0}$. 

\begin{Theorem} \StTag{DI5} Let $\CalI_2$ be a rank 3 Pfaffian system  on  a seven manifold $M$ for a second order hyperbolic PDE in the plane which is Darboux integrable with singular Pfaffian systems $\hV$ and $\cV$ satisfying $\rank \hV^\infty =2, \rank \cV^\infty = 2$.  If the Vessiot algebra for $\CalI_2$ is not ${\bf {so}}(3)$, then a local  B\"acklund transformation $\beta$ with 1-dimensional fibre between $\CalI_2$ and
$\CalI_{s=0}$ can be constructed through symmetry reduction
\begin{equation}
\beginDC{\commdiag}[3]
\Obj(0, 32)[W]{$\CalK_1 + \CalK_2$}
\Obj(0, 12)[B]{$\CalB$}
\Obj(-30, 0)[I1]{$\boxed{ \CalI_{s= 0}}$}
\Obj(30,  0)[I2]{$\boxed{ \CalI_2 }\ .$}
\mor{W}{B}{\lower 15 pt\hbox{$\bfq_{\Gamma_H}$}}[\atright, \solidarrow]
\mor{W}{I1}{$\bfq_{\Gamma_{G_1}}$}[\atright, \solidarrow]
\mor{W}{I2}{$\bfq_{\Gamma_{G_2}}$}[\atleft, \solidarrow]
\mor{B}{I2}{$\bfp_2$}[\atright, \solidarrow]
\mor{B}{I1}{$\bfp_1$}[\atleft, \solidarrow]
\enddc
\EqTag{Intro10}
\end{equation}
\end{Theorem}
The EDS $\CalK_1+\CalK_2$ in this theorem is easily constructed from $\CalI_2$. 
In the case where $F=0$ in equation \EqRef{DEFM} is a Monge-Amp\`ere equation, diagram \EqRef{Intro10} can be refined.

\begin{Theorem}\StTag{BTMA0} If $F=0$ in equation \EqRef{DEFM} is a Darboux integrable hyperbolic Monge-Amp\`ere equation, then diagram \EqRef{Intro10} can be de-prolonged producing 
\begin{equation}
\begin{gathered}
\beginDC{\commdiag}[3]
\Obj(0, 38)[KK]{$\CalK_1+ \CalK_2$}
\Obj(0, 23)[C]{$\CalB$}
\Obj(0, 8)[bB]{$\bar \CalB$}
\Obj(-35, 16)[I1]{$\CalI_1$}
\Obj(35, 16)[I2]{$\CalI_2$}
\Obj(-35, 0)[bI1]{$\bar \CalI_1$}
\Obj(35, 0)[bI2]{$\bar \CalI_2$}
\mor{KK}{I1}{$\bfq_{\Gamma_{G_1}}$}[\atright, \solidarrow]
\mor{KK}{I2}{$\bfq_{\Gamma_{G_2}}$}[\atleft, \solidarrow]
\mor{KK}{C}{$\bfq_{\Gamma_H}$}[\atright, \solidarrow]
\mor{C}{bB}{\lower 20pt \hbox{$\bfq_{\Cau_{\CalB'}}$} }[\atright, \solidarrow]
\mor{I1}{bI1}{$\bfq_{\Cau_{\CalI_{s=0}'}}$}[\atright, \solidarrow]
\mor{I2}{bI2}{$\bfq_{ \Cau_{\CalI_2'}}$}[\atleft, \solidarrow]
\mor{C}{I1}{$\bfp_1$}[\atright, \solidarrow]
\mor{C}{I2}{$\bfp_2$}[\atleft, \solidarrow]
\mor{bB}{bI1}{$\tilde \bfp_1$}[\atright, \solidarrow]
\mor{bB}{bI2}{$\tilde \bfp_2$}[\atleft, \solidarrow]
\enddc
\end{gathered} .
\EqTag{FCC}
\end{equation}
where $\bar \CalB$ is a rank 2 Pfaffian system on a 6-manifold and is a B\"acklund transformation with one dimensional fibre between the Monge-Amp\`ere representation $\bar\CalI_2$ on a 
five-dimensional manifold for the PDE $F=0$ and the Monge-Amp\`ere  system $\bar \CalI_{s=0}$ for the wave equation.
\end{Theorem}

The B\"acklund transformation $\bar \CalB$ can also be obtained using the following theorem.

\begin{Theorem}\StTag{BTMA} Let $\bar \CalI_2$ be a Monge-Amp\`ere system on a five manifold $M_2$ which is Darboux integrable (and not Monge integrable) after one prolongation and whose Vessiot algebra $\vess(\bar \CalI_2^{[1]})$ is not ${\bf so}(3)$, and let $\CalC_1+\CalC_2$ be the standard contact structure on $J^2(\real,\real)\times J^2(\real, \real)$. Then
 there exists two Lie algebras $\tilde \Gamma_{G_a}$ of infinitesimal contact transformations on $J^2(\real,\real)\times J^2(\real,\real)$ where $\tilde \Gamma_H= \tilde \Gamma_{G_1} \cap \tilde \Gamma_{G_2}$ is two 
dimensional,  which produces a (local) B\"acklund transformation $\bar \CalB$ with one dimensional fibre between $\bar \CalI_2$ and  $\bar \CalI_{s=0}$ by the reduction diagram 
\begin{equation}
\beginDC{\commdiag}[3]
\Obj(0, 32)[I]{$\CalC_1+\CalC_2$}
\Obj(0, 12)[H]{$\bar \CalB=(\CalC_1+\CalC_2) /\tilde \Gamma_H$}
\Obj(-40, 0)[I1]{$\bar \CalI_{s=0}=(\CalC_1+\CalC_2)/\tilde \Gamma_{G_1} $}
\Obj(40,  0)[I2]{$\bar \CalI_2 =(\CalC_1+\CalC_2)/\tilde \Gamma_{G_2}$\ . }
\mor{I}{H}{$\bfq_{\tilde \Gamma_H}$}[\atleft, \solidarrow]
\mor(-4,9)(-28,2){$\tilde \bfp_1$}[\atleft, \solidarrow]
\mor(4,9)(28,2){$\bfp_2$}[\atright, \solidarrow]
\mor(-2,29)(-30,5){$ \bfq_{\tilde \Gamma_{G_1}}$}[\atright, \solidarrow]
\mor(2,29)(30,5){$ \bfq_{\tilde \Gamma_{G_2}.}$}[\atleft, \solidarrow]
\enddc
\EqTag{Intro1b}
\end{equation}

\end{Theorem}

The converse of Theorem \StRef{BTMA} is proved in \cite{anderson-fels:2011a}. Given a B\"acklund transformation on a $6$-dimensional
manifold between a Monge-Amp\`ere system $\bar\CalI_2$ which is Darboux integrable (and not Monge integrable) after prolongation, and the Monge-Amp\`ere system for the wave equation, there exists a three dimensional Lie algebras of vector-fields $\Gamma_{G_2}$, and  two-dimensional sub-algebra $\Gamma_H \subset \Gamma_{G_2} $ such that the B\"acklund transformation $ \CalB$ can be constructed by symmetry reduction. In particular,  the equation
\begin{equation}
	u_{xy}= \frac{ \sqrt{1-u_x^2} \sqrt{1-u_y^2} }{\sin u },
\EqTag{VGso3}
\end{equation}
has Vessiot algebra ${\bf so}(3)$ which has no two-dimensional subalgebra. Therefore equation \EqRef{VGso3} admits no B\"acklund transformation with one-dimensional fibre with the wave equation. This disagrees with Theorem 1 in \cite{Clelland-Ivey:2009a}. The article \cite{Clelland-Ivey:2009a} also establishes the existence of B\"acklund  transformations through the Cartan-K\"ahler theorem while the B\"acklund transformations whose existence is given here require only $C^\infty$ group actions.



\section{Reduction of Differential Systems}

\subsection{Preliminaries}

An  EDS $\CalI \subset \Omega^*(M)$ on a manifold $M$ is a differentially closed ideal (see \cite{bryant-chern-gardner-griffiths-goldschmidt:1991a}, and \cite{ivey-langsberg:2003}).  The EDS $\CalI$ has constant rank if each  of its homogeneous components $\CalI^k \subset \Omega^k(M)$  coincide with the sections 
	$\CalS(I^k)$ of a  constant rank sub-bundle $I^k \subset \Lambda^k(T^*M)$. If $\CalA$ is a subset of $\Omega^*(M)$, 
	we let  $\langle \CalA  \rangle_{\text{alg}} $ and $\langle \CalA \rangle_{\text{diff}} $ be the algebraic and differential ideals  generated by $\CalA$. 
		
A constant rank Pfaffian differential system $\CalI$ is an EDS for which there exists  a sub-bundle $I\subset T^*M$  such that $\CalI = \langle \CalS(I) \rangle_\text{diff}$.  We also refer to a constant rank sub-bundle $I\subset T^*M$ as a Pfaffian system.

A {\deffont local first integral}  of a Pfaffian system $I$ is smooth function 
	$f\: U \to \Real$, defined on an open set $U$,  and such that $d f \in I$. For each point $x \in M$ we define
\begin{equation}
	I_x^\infty  = \{\, df _x \, | \,   \text{$f$ is a local first integral, defined about $x$}\,\} .
	\EqTag{Iinfty}
\end{equation}
We shall always assume that  $I^\infty   = \cup _{x\in M} I_x^\infty$ is a constant rank bundle on $M$.
It is easy to verify check that $I^\infty$ is  the  (unique) maximal,  completely integrable,  Pfaffian subsystem of  $I$.
	The bundle $I^\infty$ can be computed algorithmically from the derived  flag of $I$.  The derived system $I'\subset I$ of a Pfaffian system $I$ is defined pointwise by
$$
I'(p)= span \{ \theta_p \ | \  \theta \in \CalS(I) \quad d \theta \equiv 0 \quad \mod \ I \}.
$$
The system $I$ is integrable if it satisfies the Frobenius condition if $I'=I$. Letting $I^{(0)}=I$ and assuming $I^{(k)}$ is constant rank we define inductively
$$
I^{(k+1)} = (I^{(k)})' , \quad k=0,1\ldots, N
$$
where $N$ is the smallest integer where $I^{(N+1)} = I^{(N)}$. Therefore  $I^\infty = I^{(N)}$ when $I^{(k)}$ are constant rank.
More information about this sequence can be found in \cite{bryant-chern-gardner-griffiths-goldschmidt:1991a}, and \cite{ivey-langsberg:2003}.

We also recall the definition of an integrable extension \cite{bryant-griffiths:1995a}. First,
if $\CalI$ and $\CalJ$ are differential systems, we let  $\CalI +\CalJ =  \langle \CalI \cup \CalJ \rangle_{\text{alg}}$ which is also a differential system. Let $\bfp :M \to N$ be a surjective submersion and $\CalI$ an EDS on $N$.  An EDS  $\CalE$ on $M$
	is called an integral extension of  $\CalI$ if there exists a subbundle  $J \subset  \Lambda^1(M)$ of rank $\dim M - \dim N$, such that $J$ is transverse to  $\bfp$, that is,
\begin{subequations}
\EqTag{IntExt2}
\begin{equation}
\ann(J) \cap \text{ker } (\bfp_*)=0 \, ,
\EqTag{IntExt2a}
\end{equation}
and
\begin{equation}
	\CalE = \langle\,  \bfp^*(\CalI)  +  \CalS(J) \rangle_{\text{alg}}.
\EqTag{IntExt2b}
\end{equation}
\end{subequations}
	A sub-bundle $J$  satisfying these two properties is called an {\deffont admissible } sub-bundle for the extension $\CalE$.  Given any (immersed) integral manifold $s:P \to M$ of $\CalE$, then by condition \EqRef{IntExt2a} the composition $\bfp \circ s:P \to N$ is an (immersed) integral manifold of $\CalI$.
If  $s:P \to N$ is an integral manifold of $\CalI$, then by condition \EqRef{IntExt2b} the restricted of the EDS $\CalE$ to $\bfp^{-1}(s(P))$ is an integrable Pfaffian system.



\subsection{Reduction, Pullback Bundles, and Semi-basic Forms} \StTag{prelim}

Let $\bfp:M\to N$ be a smooth submersion and let $\CalI$ be an EDS on $M$ and $\CalI/\bfp$ be the reduction of $\CalI$ (see \EqRef{SQ}). In this section we give easily verifiable conditions which guarantee that $\CalI/\bfp$ is constant rank and provide a simple way to compute local basis of sections for $\CalI/\bfp$.

Let $Vert(M) = \ker (\bfp_*) \subset T^*M$ be the vertical distribution for the submersion $\bfp:M\to N$. A vector field $X$ on $M$ taking values in $Vert(M)$ is called a vertical vector-field. A differential form $\theta\in \Omega^*(M) $ is $\bfp$ {\bf semi-basic} if $X \hook \theta= 0 $ for all vertical vector-fields $X$, and $\theta$ is called $\bfp$ {\bf basic} if there exists $\bar \theta \in \Omega^*(N)$ such that $\theta = \bfp^*(\bar \theta)$. 
Likewise, for any sub-bundle $I^k \subset \Lambda^k(T^*M)$ the subset of $\bfp$ semi-basic $k$-forms is defined by
\begin{equation}
I^k_{\bfp,sb} = \{ \ \theta \in  I^k \ | \ X \hook \theta =0 \ {\text for \ all } \  X \in Vert(M) \}.
\EqTag{Ibsb}
\end{equation}
Clearly any $\bfp$ semi-basic $k$-form takes values in   $I^k_{\bfp,sb} $.




Now let $\barI^k \subset \Lambda^k(N)$ be a rank $r$ sub-bundle. The pullback bundle $\bfp^*(\barI^k) \subset \Lambda^k(T^*M)$ is the rank $r$ sub-bundle given point-wise by
\begin{equation}
\bfp^*(\barI^k)_x = \{\ \bfp^* (\bar \theta_y) \ | \ {\text for \ all }  \ \bar \theta_y \in  \Lambda^k(T^*_yN)\ , \ y=\bfp(x) \  \}.
\end{equation}
Note that $\bfp^* (\barI^k) \subset \Lambda^k(T^*M)_{\bfp,sb}$.  Lemma \StRef{PB2} below characterizes pullback bundles and is the bundle version of Proposition 6.1.19 in \cite{ivey-langsberg:2003}.

\begin{Lemma} \StTag{PB2} Suppose the fibres of $\bfp:M \to N$ are connected. Given a sub-bundle $I^k \subset  \Lambda^k(T^*M)$ there exists a sub-bundle $\barI^k \subset \Lambda^k(T^*N)$ such that $I^k = \bfp^*(\barI^k)$ if and only if  $I^k$ is $\bfp$ semi-basic and for all
vertical vector fields $X$
\begin{equation}
\Psi^{X,*}_{t}( I^k)=  I^k\, ,
\EqTag{Invs0}
\end{equation}
where $\Psi^X_t$ is the flow of $X$. \end{Lemma}

Condition \EqRef{Invs0} is equivalent to the condition that for all $\theta \in \CalS(I^k) $ 
\begin{equation}
X \hook \theta =0 \ ,\quad L_X \theta\in \CalS(I^k) \ ,
\EqTag{Invs1}
\end{equation}
for all vertical vector-fields $X$. See Proposition 6.1.19 in \cite{ivey-langsberg:2003}. Note that $\Lambda^k(T^*M)_{\bfp,sb}$
satisfies \EqRef{Invs0} (or \EqRef{Invs1}) and that $\Lambda^k(T^*M)_{\bfp,sb}= \bfp^*( \Lambda^k(T^*N))$.

In analogy with the definition of $\CalI/\bfp$ in \EqRef{SQ} the reduction of a bundle $I^k\subset \Lambda^k(T^*M)$ is defined by
\begin{equation}
I^k/\bfp = \{ \theta \in \Lambda^k(T^*N) \ | \ \bfp^*(\theta ) \in I^k \}\, .
\EqTag{Ibmodp}
\end{equation}
Lemma \StRef{PB2} then has the following corollary.

\begin{Corollary} \StTag{CR1} The set $I^k/\bfp$ is a constant rank bundle if and only if $I^k_{\bfp,sb}$ (equation \EqRef{Ibsb})  is constant rank and vertically invariant (equation \EqRef{Invs0}). In which case $I^k_{\bfp,sb}=\bfp^* (I^k/\bfp)$.
\end{Corollary}



Let $\CalI$ be a constant rank EDS. The subset $\CalI_{\bfp,sb} \subset \CalI $ of semi-basic forms is
\begin{equation}
\CalI_{\bfp,sb} = \{ \ \theta \in \CalI \ | \ X \hook \theta= 0 \ {\text for \ all \ vertical \ vector \ fields }\ X \ \}.
\EqTag{Isb}
\end{equation}
If $\CalI/\bfp$ is the reduction of $\CalI$ then
$$
\bfp^*(\CalI/\bfp) \subset \CalI_{\bfp,sb}.
$$
Note that  neither $\bfp^*(\CalI/\bfp)$ or $\CalI_{\bfp,sb}$ are differential ideal. 
The next theorem is the basic theorem on reduction. It's proof follows immediately from Lemma \StRef{PB2} and Corollary \StRef{CR1}.

\begin{Theorem}\StTag{RIp} Let $\bfp:M \to N$ be a smooth submersion with connected fibres. If $\CalI$ is a constant rank EDS and if $I^k_{\bfp,sb}$ (see equation \EqRef{Ibsb}) are constant rank and satisfy the vertical invariance condition in equation \EqRef{Invs0} (or \EqRef{Invs1})  then

\noindent
{\bf [i]} there exists  constant rank bundles $\barI^k \subset \Lambda^k(T^*N)$ such that  $ I^k_{\bfp,sb} = \bfp^* (\barI^k) $,

\noindent
{\bf [ii]} $\CalI/\bfp$ is a constant rank EDS (with bundles $\barI^k$), and

\noindent
{\bf [iii]} $\langle \CalI_{\bfp,sb} \rangle_{\text alg} = \langle \bfp^* (\CalI/\bfp) \rangle_{\text alg}$ .
\end{Theorem}

The vertical invariance conditions in Theorem \StRef{RIp} can be inferred from the invariance of $\CalI$.

\begin{Lemma}\StTag{IndInv} If $\CalI $ is constant rank and vertically invariant (equation \EqRef{Invs0}) then

\noindent
{\bf [i]} the semi-basic subset  $I^k_{\bfp,sb}$ is vertically invariant for each $k$,

\noindent
{\bf [ii]} the space of invariants $I^\infty$ is vertically invariant and 

\noindent
{\bf [iii]} the derived systems $ I^{(k)}$  are vertically invariant.
\end{Lemma}

\begin{proof} Let $X$ and $Y$ be vertical vector fields on $M$, and let $\Psi^X$ be the flow of $X$. Then
$$
Y\hook (\Psi^X_t)^* \theta = \theta( (\Psi^X_{t,*}) Y ) .
$$
Now $\Psi^X_{t,*} Y$ is a vertical vector-field (on its domain), and so if $ \theta \in I^k_{\bfp,sb}$ then $ \theta( (\Psi^X_{t,*}) Y )  =0 $
which proves part {\bf [i]}. 

Both  {\bf [ii]} and the case if $I^{(1)}=I'$ in part {\bf [iii]} follow from the vertical invariance of $I$ (equation \EqRef{Invs0}), and the
commutativity of the exterior derivative and the pullback. The rest of part {\bf [iii]} follows by induction
\end{proof}

The following lemma is particularly useful for computing local basis of sections in reduction.

\begin{Lemma}\StTag{PB3} Suppose $I^k_{\bfp,sb} = \bfp^*(\barI^k) $, and let $\sigma :\barU \to U$ be a local cross-section to $\bfp$. Then

\noindent
{\bf [i]} $ \sigma^*(I_{\bfp,sb}^k) = \barI^k$, and

\noindent
{\bf [ii]}  if $\theta \in \CalS(I^k_{\bfp,sb}|_U)$ then  $\sigma ^* \theta \in \CalS(\barI^k|_{\barU})$.

\noindent
{\bf [iii]} Given a local basis of sections $\{ \theta^a \}$ for $I^k_{\bfp, sb}|_U$ and a cross-section $\sigma :\barU \to U$, then $\sigma^* (\theta^a)$ is a local
basis of sections for $\barI^k|_{\barU} $. 

\noindent
{\bf [iv]} About each point $x\in M$ there exists an open set $U$ and a local basis $\{\theta^a\}$ of sections for $I^k_{\bfp,sb}$ where each $\theta^a$ is $\bfp$ basic.
\end{Lemma}

\begin{Remark} The point behind Lemma \StRef{PB3} is the following. Given generators $\CalI= \langle \alpha^1,\ldots, \alpha^k \rangle_{\text{alg}}$ we find generators for the semi-basic forms $\CalI_{\bfp,sb}=\spn \{ \tilde \alpha^1,\ldots, \tilde \alpha^r \}$   from $\CalI$ using only linear algebra (equation \EqRef{Isb}). Local generators for $\CalI/\bfp$ are then computed using the pullback $(\CalI/\bfp)_{\barU} = \langle \sigma^* \tilde \alpha^1, \ldots, \sigma ^* \tilde \alpha^ r \rangle_{\text{alg}}$ by a local section $\sigma: \barU \to U$ of $\bfp:M \to N$. Lemma \StRef{PB3} is also useful for producing an adapted set of generators for $\CalI$. In particular  if
$(\CalI / \bfp )_{\barU} = \langle  \bar \alpha^1, \ldots,  \bar \alpha^ r \rangle_{\text{alg}}$ then
$$
\CalI|_U = \langle \bfp^* \bar \alpha^1, \ldots,  \bfp^* \bar \alpha^ r , \hat \alpha^{r+1} ,\ldots, \hat \alpha^ k \rangle_{\text{alg}},
$$
where the forms $ \bfp^* \bar \alpha^s $ are $\bfp$ basic.
\end{Remark}

\subsection{Reduction by Cauchy Characteristics} \StTag{S23}

Here we briefly review the theory of reduction by Cauchy characteristics again see \cite{bryant-chern-gardner-griffiths-goldschmidt:1991a}, \cite{ivey-langsberg:2003} and relate this to the previous section.  

Let ${\Delta} \subset TM$ be an integrable distribution. Given $x\in M$ we denote by $S_x$ the maximal (connected) integral manifold through $x$.  Let $\bfq_\Delta:M \to M/\Delta$ be the quotient map by the maximal integral manifolds of $\Delta$ and so $\bfq_\Delta(x) = S_x$.
We call ${\Delta}$ {\deffont regular} if $M/\Delta$, endowed with the quotient topology, admits a unique manifold structure such that $\bfq_\Delta$ is a smooth submersion.  Note that regularity implies $\Delta$ is constant rank.

Let $\CalI$ be a constant rank EDS. The set  
\begin{equation}
\Cau_{\CalI} = \Cau(\CalI) = \{ X  \in TM \ | \ X \hook \CalI \subset \CalI \} 
\EqTag{CauD}
\end{equation}
is the space of Cauchy characteristic of $\CalI$ which is an integrable distribution \cite{bryant-chern-gardner-griffiths-goldschmidt:1991a}. From here on we will assume that $A(\CalI)$ is regular and denote by $\bfq_{A}:M \to M/A(\CalI)$  the smooth quotient map to the leaf space.

The main theorem on the reduction by Cauchy Characteristics is the following, see \cite{bryant-chern-gardner-griffiths-goldschmidt:1991a}, \cite{ivey-langsberg:2003}.

\begin{Theorem} \StTag{CQ} Let $\CalI$ be a constant rank EDS on an $n$-dimensional manifold $M$ and suppose $A(\CalI)$ is regular.  If  $I^k_{\bfq_A, sb }\subset I^k$ is constant rank for each $k=1,\ldots, n$,  then the constant rank EDS  $\bar \CalI =\CalI/\bfq_\Cau $ statisfies
$$
\CalI = \langle \bfq_A^* \bar \CalI \rangle_{alg},
$$
\end{Theorem}

Equation \EqRef{CauD} implies that $\CalI$ is vertically invariant (equation \EqRef{Invs0}), and therefore by Lemma \StRef{IndInv}, the hypothesis in Theorem \StRef{RIp} are satisfied and $\CalI/\bfq_A$ is a constant rank EDS. The difference between Theorem \StRef{CQ} and \StRef{RIp}, which is an essential point in reduction by Cauchy characteristics, is that $\CalI$ admits a set of $\bfq_A$ semi-basic generators.

 Again, local generators for $\CalI/\bfq_A$ can be constructed using Lemma \StRef{PB3}. For examples of Cauchy characteristic reduction see \cite{bryant-chern-gardner-griffiths-goldschmidt:1991a},  \cite{ivey-langsberg:2003}.

\begin{Remark} \StTag{CDelta} First note that Theorem \StRef{CQ} holds for any regular integrable distribution $\Delta \subset A(\CalI)$, where $\bfq_\Delta:M \to M/\Delta$ is the smooth quotient map. Let $\bfp :M \to N$ be a smooth submersion with connected fibres and $\CalI$ and EDS on $M$. If the EDS $\langle \CalI_{\bfp,sb} \rangle_{\text{alg}}$ is vertically invariant (equations \EqRef{Invs0}, \EqRef{Invs1}) then 
$$
Vert(M) \subset A( \langle \CalI_{\bfp,sb} \rangle_{\text{alg}} )
$$
is an integrable sub-distribution of the Cauchy characteristics of $\langle \CalI_{\bfp,sb} \rangle_{\text{alg}}$. Then 
$$
\CalI/\bfp = \left( \langle \CalI_{\bfp,sb} \rangle_{\text{alg}} \right) /\bfp=  \left( \langle \CalI_{\bfp,sb} \rangle_{\text{alg}} \right) /Vert(M) .
$$
In other words the reduction $\CalI/\bfp$ is Cauchy characteristic reduction of the EDS generated by the semi-basic forms by the subset $Vert(M)$. 
\end{Remark}

\subsection{Group Reduction}

Let  $G$ be a Lie group acting on $M$. The action $\mu:G\times M \to M$ is said to be regular if the space of orbits $M/G$ with the quotient topology admits a smooth manifold structure such that the projection map $\bfq_G : M \to M/G$  is a submersion.  We let $\Gamma_G$ denote the Lie algebra of infinitesimal generators of $G$ and let $\bfGamma_G$ be the distribution spanned by $\Gamma_G$. In particular for regular actions $Vert(M) = \ker (\bfq_{G,*}) = \bfGamma_G$.  The action of $G$ is a symmetry group of the EDS $\CalI$ if $g^* \CalI = \CalI$ for all $g\in G$. 

We assume the action of $G$ on $M$ is regular and a symmetry group of $\CalI$ and we shall write $\CalI/G$ in place of $\CalI/\bfq_G$. Equation \EqRef{SQ} gives
\begin{equation}
\CalI/G = \{ \ \bar \theta \in \Omega^*(M/G) \ | \ \bfq_G^* \bar \theta \in \CalI \ \} \, .
\end{equation}

The analogue to Corollary 2.3 in \cite{bryant-chern-gardner-griffiths-goldschmidt:1991a}, which characterizes basic forms, and Lemma \StRef{PB2}  is the following.

\begin{Lemma}  A differential form $\theta \in \Omega^k(M)$ is $\bfq_G$ basic if and only if 
$$
X \hook \theta = 0, \quad  g^* \theta= \theta, 
$$
for all $X\in \bfGamma$. A rank $r$ sub-bundle $I^k \subset \Lambda^k(M)$
satisfies 
$$
X\hook I^k = 0 , \ g^* I^k=I^k 
$$
if and only if there exists a rank $r$ sub-bundle $\barI^k \subset \Lambda^k(M/G)$ such that $I^k = \bfq_G^*(\barI^k)$.
\end{Lemma}

The set  $\CalI_{G,sb} \subset \CalI$ of $\bfq_G$ semi-basic forms is
$$
\CalI_{\bfq_G,sb} = \{ \ \theta \in \CalI \ | \ X \hook \theta =0 \ \}.
$$ 
It follows in the same manner as in Lemma \StRef{PB3} that the subset $\CalI_{\bfq_G,sb}$ is $G$-invariant. The analogue of Theorem \StRef{RIp} for symmetry groups is the following.

\begin{Theorem} \StTag{RG} Let $G$ be a Lie group acting regularly on $M$ which is a symmetry group of the constant rank EDS  $\CalI$.  If the bundles $I^k_{\bfq_G,sb} $ are constant rank, then

\noindent
{\bf [i]} there exists  constant rank bundles $\barI^k \subset \Lambda^k(T^*(M/G))$ such that  $ I^k_{\bfq_G,sb} = \bfq_G^* (\barI^k) $,

\noindent
{\bf [ii]} $\CalI/\bfq_G$ is a constant rank EDS (with bundles $\barI^k$), and

\noindent
{\bf [iii]} $\langle \CalI_{\bfq_G,sb} \rangle_{\text alg} = \langle \bfq_G^* (\CalI/\bfq_G) \rangle_{\text alg}$ .

\end{Theorem}
The vertical invariance condition \EqRef{Invs0} is not necessary in Theorem \StRef{RG} because of the hypothesis that $G$ is a symmetry group of $\CalI$. Also note, as in Remark \StRef{CDelta}, that the vertical space $\bfGamma_G $ is a subset of the Cauchy characteristics for $\langle \CalI_{\bfq_G,sb} \rangle_{alg}$. If $G$ has connected orbits then $\CalI/G$ is Cauchy reduction of $\langle \CalI_{\bfq_G,sb} \rangle_{\alg}$ by $\bfGamma_G$. 




A symmetry group $G$ of an EDS $\CalI$ is said to be {\deffont transverse} to $\CalI$ if
\begin{equation}
	 \ann(I^1) \cap \bfGamma_G = {0},
\EqTag{Itrans}
\end{equation}
where $I^1\subset T^*M$ is the sub-bundle whose sections are the differential 1-forms in $\CalI$. The actions considered in this article will all satisfy the transversality condition \EqRef{Itrans}. For more information on transversality in reduction see  \cite{anderson-fels:2005a}. In particular in \cite{anderson-fels:2005a} it is shown that for transverse actions the hypothesis in Theorem \StRef{RG} are satisfied and $\CalI/G$ is a constant rank EDS. For example, by equation \EqRef{Itrans} the rank of the one forms in $\CalI/G$ is given by $ \rank\, ( \CalI/G)^1 = \rank\, \CalI^1 - \rank \bfGamma_G$.

\subsection{The Derived System} \StTag{S22}

We now examine the behavior of the the derived bundles for a Pfaffian system $I \subset T^*M$, under reduction.

\begin{Theorem}\StTag{DRED} Let $I \subset T^*M$ be a constant rank sub-bundle and assume that $\barI = I/\bfp \subset T^*N$  is constant rank. Then for any $x\in M$,  $I'_{\bfp,sb}(x) = \bfp^* (\barI'(y))$ where $y =\bfp(x)$, and
\begin{equation}
\rank \barI '(y)  \leq \rank I'(x) \leq \rank \barI'(x)+ \rank I  - \rank \barI  
\EqTag{DRKS}
\end{equation}
\end{Theorem}

\begin{proof} Let $x\in M$, $y=\bfp(x)$ and let $\barU$ be an open set in $N$ containing $y$ where  $\bar \theta^i,\bar \eta^a$ form a local basis of sections for $\barI $ such that
\begin{equation}
\barI|_{\barU} = span \{ \bar \theta^i, \bar \eta^a \}, \quad \barI'(q) =span \{ \bar \eta^a_q \} .
\EqTag{CF1}
\end{equation}
Next let $U\subset \bfp^{-1}(\barU)$ be an open subset containing $x$ where $\bfp^* \bar\theta^i, \bfp^*\bar \eta^a$ can be completed to
a local basis of sections for $I$ by $\theta^ r$,
$$
I|_U = span \{ \theta^r, \bfp^* \bar\theta^i , \bfp^* \bar \eta^ a\} .
$$

If $\rho \in I'_{\bfp,sb}(x)$, then there exists $A_r,B_i,C_a \in C^\infty(U)$ such that
\begin{equation}
\rho = A_r\, \theta^r + B _i\, \bfp^*\bar \theta ^i + C_ a \, \bfp^*\bar \eta^a
\EqTag{R1}
\end{equation}
and $A_r(x) =0$ by the condition $\rho$ is $\bfp$ semi-basic at $x$.  Taking into account that $A_r(x)=0$, the condition $\rho \in I'(x)$ is
\begin{equation}
B_i(x)\, \bfp^* d \bar \theta^i_x +C_a(x) \, \bfp^* d \bar \eta^a_x =0 \quad \mod I .
\EqTag{drho2}
\end{equation}
Equation \EqRef{drho2} implies there exists $\alpha_r,\beta_i,\gamma_a \in T_x^*M$ so that
\begin{equation}
B_i(x)\, \bfp^* d \bar \theta^i_x + C_a(x)\,  \bfp^*d \bar \eta^a_x = \alpha_r \wedge \theta^r_x+ \beta_i \wedge \bfp^* \bar \theta^i_x + \gamma_a \wedge \bfp^* \bar \eta^a_x
\EqTag{drho3}
\end{equation}
where we may assume $\alpha_r$ contain no $\bfp^* \bar \theta^i_x $ or $\bfp^* \bar \eta^a_x$ terms.

Now $B_i(x) \, \bfp^* d \bar \theta^i_x+ C_a(x)\,  \bfp^* d\bar \eta^a_x $ is $\bfp$ semi-basic, and therefore since $\theta^r_x$ are now $\bfp$ semi-basic the term $\alpha_r \wedge \theta^r_x$ is zero. Furthermore the terms $\beta_i,\gamma_a$ are $\bfp$ semi-basic which implies there exists $\bar \beta_i, \bar \gamma_a \in T^*_yN$ so that $\beta^i=\bfp^* \bar \beta^i$, $\gamma^a = \bfp^* \bar \gamma^a$.  Equation \EqRef{drho3} then reads 
$$
B_i(x) \, d \bar \theta^i_y + C_a(x) \, d  \bar \eta^a_y = \bar \beta_i \wedge  \bar \theta^i_y + \bar \gamma_a \wedge \bar \eta^a_y .
$$
This equation implies $ B_i(x) \, d \bar \theta^i_y  \in \bar I'|_y$ which by the choice of coframe for $\barI$ in equation \EqRef{CF1}  gives $B_i(x) = 0 $. Therefore $\rho$ in equation \EqRef{R1} satisfies $\rho(x) =C_a(x)\,  \bfp^* \bar \eta^a_y$. This implies $I'_{\bfp,sb}(x)=\bfp^*\barI'(y) $. 

Finally to show the  inequality \EqRef{DRKS} note that the index $r$ in equation \EqRef{R1}  satisfies $1\leq r \leq \rank I -\rank\bar I $. This shows the  inequality \EqRef{DRKS} holds. 
\end{proof}

\begin{Corollary}\StTag{DREDC} The subset $I'_{\bfp,sb}$ is constant rank if and only if $\barI'$ is constant rank in which case the ranks are equal, $I'/\bfp =\barI'$, and
$$
0 \leq \rank I' - \rank \barI ' \leq \rank I - \rank \barI .
$$
\end{Corollary}

Theorem \StRef{DRED} and Corollary \StRef{DREDC} state that bundle reduction and derivation commute. We also have the following corollary of Theorem \StRef{DRED}.

\begin{Corollary} \StTag{RDone} Let $\CalI=\langle \CalS(I)\rangle_{\text diff}$ be the Pfaffian system generated by $I\subset T^*M$, and suppose $\bar\CalI= \langle \CalS(\barI)\rangle_{\text diff}$ is the Pfaffian system generated by $\barI=I/\bfp$. Then 
\begin{equation}
\CalI'/\bfp \cap \Omega^1(N) = \bar \CalI' \cap \Omega^1(N) , \quad{\text and} \quad \bar \CalI' \subset \CalI'/\bfp .
\EqTag{Qdone}
\end{equation}
\end{Corollary}

\begin{Remark} In Corollary \StRef{RDone}, by definition $\bar \CalI'$ is the Pfaffian system generated by $\barI '$, while $\CalI'/\bfp$ is
not necessarily a Pfaffian system and so there is no reason for these to be equal. Corollary \StRef{RDone} does say that these two systems have the same set of one-forms. 
\end{Remark}

Using definition \EqRef{Iinfty} we now look at the reduction of $I^\infty$. We begin with the following lemma.

\begin{Lemma}\StTag{CIred} Suppose that $I\subset T^*M$ is a constant rank completely integrable sub-bundle, and that $\barI = I/\bfp$ is constant rank. Then $\barI$ is completely integrable.
\end{Lemma}
\begin{proof}  Since $I$ is completely integrable $I'=I$ and so $I'_{\bfp,sb}$ is constant rank and $I'_{\bfp,sb} = \bfp^* \barI'$. 
However $\barI'=I'/\bfp = I/\bfp=\barI $ and so $\barI$ is completely integrable.
\end{proof}

\begin{Corollary}\StTag{Rint} Suppose that $I$, $I^\infty$,  $\barI = I/\bfp$, and $I^\infty_{\bfp,sb}$ are constant rank.  Then
$$
I^\infty/\bfp = \barI ^\infty.
$$
\end{Corollary}
\begin{proof}   We begin by using Lemma \StRef{CIred}  to note that  $I^\infty/\bfp$ is completely integrable. Now 
$\bfp^* \barI^\infty \subset I^\infty$ implies 
\begin{equation}
\barI^\infty \subset I^\infty/ \bfp.
\EqTag{Rinf2}
\end{equation}
However $I^\infty/\bfp \subset I/\bfp$ which by the maximality of $\barI^\infty$ implies $I^\infty/\bfp \subset \barI^\infty$. This
along with equation \EqRef{Rinf2} imply $\barI^\infty = I^\infty/\bfp$. 
\end{proof}

The derived series behaves similarly. 

\begin{Corollary} Let $I^k$ be the  $k^{th}$ derived bundle and suppose  $I^{k+1}$ is constant rank, $I^k/\bfp$ is constant rank, and $(I^k/\bfp) '$ is
constant rank. Then $I^{k+1}/\bfp$ is constant rank and $(I^k/\bfp)' = I^{k+1}/\bfp$. Furthermore, if $I^k/\bfp$ are constant rank for $k=0\ldots m$, and $(I/\bfp)^k$ are constant rank for $k=0\ldots m$, then $I^{k}/\bfp$ is constant rank and $(I/\bfp)^{k}=I^k/\bfp$, $k=0\ldots m$.
\end{Corollary}

\section{B\"acklund Transformations} 

We recall the following key result proved in \cite{anderson-fels:2011a}, from which Theorem \StRef{ThIntro1} follows.

\begin{Theorem}\StTag{ORD} Let $G$ and $H$ be Lie groups acting act regularly on $M$ where $H$ is a subgroup of $G$. If $G$ is a symmetry group of the EDS $\CalI$ which acts transversally to $\CalI$, then 
$$
\begin{gathered}
\beginDC{\commdiag}[3]
\Obj(-25, 0)[I]{$(M,\CalI)$}
\Obj(0, 0)[H]{$(M /H, \CalI/H)$}
\Obj(0, -17)[I1]{$(M /G, \CalI/G)$}
\mor{I}{H}{$\bfq_H$}[\atleft, \solidarrow]
\mor{H}{I1}{$\bfp$}[\atleft, \solidarrow]
\mor{I}{I1}{$\bfq_{G}$}[\atright, \solidarrow]
\enddc
\end{gathered}\ ,
$$
is a commutative diagram of integrable extensions where
\begin{equation}
\bfp : M/H \to  M/G \quad\text{is defined by} \quad\bfp(Hx) = Gx .
\EqTag{IntExt14}
\end{equation}
\end{Theorem}

Applying Theorem \StRef{ORD} with the hypothesis of Theorem \StRef{ThIntro1} proves Theorem \StRef{ThIntro1} and  provides the B\"acklund transformation $\CalI/H$ in diagram \EqRef{Intro1}. We demonstrate Theorem \StRef{ORD} and Theorem \StRef{ThIntro1} with an example.

\begin{Example} \StTag{EX1} Let $\CalI = \CalK_1+\CalK_2$ be the standard contact system on $J^3(\real,\real)\times J^3(\real,\real)$. In local coordinates $(x,v,v_x, v_{xx},v_{xxx}; y,w, w_y, w_{yy}, w_{yyy})$ we have
\begin{equation}
\begin{aligned}
\CalI = \langle\  &&  \theta^v& = dv -v_x dx,\  &\theta^v_x & = d v_x -v_{xx} dx,\ & \theta^v_{xx}&= dv_{xx}-v_{xxx}dx, & \ \\
 & & \theta^w&=dw -w_ydy,\ &\theta^w_y&= d w_y -w_{yy} dy,\ & \theta^w_{yy}&=dw_{yy}-w_{yyy}dy &\  \rangle_{\text diff} .
\end{aligned}
\EqTag{Twof}
\end{equation}
Let $\Gamma_{G_2}$ be the prolongation of Lie algebra of infinitesimal generators for the local action of $SL(2)$ acting diagonally on $v$ and $w$,
\begin{equation}
\Gamma_{G_2}= {\text span} \{\ \partial_v - \partial_w, \ {\text pr}(v \partial_v+w \partial_w), \ , {\text pr}( v^2\partial_v - w^2\partial_w)\ \} ,
\EqTag{GG2}
\end{equation}
and let
\begin{equation}
\begin{aligned}
\Gamma_{G_1} &= {\text span} \{ \partial_v, \partial_w, pr(v\partial_v + w\partial_w) \} , \\
\Gamma_H &= \Gamma_{G_1}\cap \Gamma_{G_2}=  {\text span} \{ \partial_v- \partial_w, pr( v \partial_v +w \partial_w)  \}.
\end{aligned}
\EqTag{GG1H}
\end{equation}
Let $M$ be the open set $v,v_x,w,w_y> 0$.  The distributions $\Gamma_H$, and  $\Gamma_{G_a}$ are regular on $M$. The manifolds $M/\Gamma_{G_a}$ are 7 dimensional while $M/\Gamma_H$ is 8 . With the following choices of coordinates
\begin{equation}
\begin{aligned}
M/\Gamma_{G_1} &=  \ (x,y,z,z_x,z_y,z_{xx},z_{yy}) ,  \quad
M/\Gamma_{G_2}  =  \  (x,y,u,u_x,u_y,u_{xx},u_{yy}) , \ \\
M/\Gamma_H  &=  \ (x,y,V,W,V_x,W_y,V_{xx},W_{yy}),
\end{aligned}
\EqTag{COORD1}
\end{equation}
the quotient maps $\bfq_{\Gamma_{G_1}}$, $\bfq_{\Gamma_{G_2}}$ and $\bfq_{\Gamma_H}$ in coordinates \EqRef{COORD1} are
\begin{equation}
\begin{aligned}
& \bfq_{\Gamma_{G_1}}= (x=x,=y,z=\log(\frac{v_x}{w_y}), z_x=D_x(z)=\frac{v_{xx}}{v_x}, z_y=-D_y(z)=\frac{w_{yy}}{w_{y}}, z_{xx}=D_x(z_x), z_{yy}=D_y(z_y) ), \\
&\bfq_{\Gamma_{G_2}}= (x=x,y=y, u  =\log  \frac{2w_y v_x}{(v+w)^2}\ , u_x  =D_x(u) , \
u_y  = D_y(u) , u_{xx}=D_x(u_x), u_{yy}=D_y(u_y) ) ,\\
& \bfq_{\Gamma_{H}}  =\  (x=x,y=y,V=\log \frac{v_x}{v+w},W=\log \frac{w_y}{v+w}, V_x = D_x(V), W_y = D_y(W),\qquad \qquad \qquad \\ 
& \qquad \qquad  V_{xx}=D_x(V_x),W_{yy}=D_y(W_y)),
\end{aligned}
\EqTag{bfp1}
\end{equation}
where $D_x$ and $D_y$ are the total derivatives.  
Using the maps in equation \EqRef{bfp1}, the projection maps $\bfp_1 $ and $\bfp_2$ which make the commutative diagram \EqRef{Intro10}  are easily
seen to be given in coordinates by
\begin{equation}
\begin{aligned} 
&\bfp_1=(x=x,y=y, z= V-W,  z_x = V_x+e^V,z_y = -W_y-e^W, z_{xx}  = D_x(z_x), z_{yy} = D_y(W_y) ),\\
&\bfp_2=(x=x,y=y, u=\log (2)+V+W, u_x = V_x- e^V , u_y= W_y- e^W, u_{xx}=D_x(V_x), u_{yy}=D_y(W_{y})). 
\end{aligned}
\EqTag{DP12}
\end{equation}

Since $\Gamma_{G_a}$ and $\Gamma_H$ consist of prolongations of vector-fields, they are infinitesimal symmetries of $\CalI$ and so we can
use Lemma \StRef{PB3} to compute the reductions that produce diagram \EqRef{Intro1}
\begin{equation}
\begin{aligned}
\CalI/\Gamma_{G_1} & = \langle du - z_x dx-z_y dy, dz_x -z_{xx} dx, dz_y-z_{yy} dy\rangle_{\text diff}   \\
\CalI/\Gamma_{G_2} & =\langle du - u_x dx-u_y dy, du_{x} -u_{xx}dx -e^u dy, du_y-e^u dx-u_{yy}dy\rangle_{\text diff} \\
\CalI/\Gamma_H& =\langle  dV-V_x dx+e^W dy,
dW+e^Vdx-W_ydy,  \\
& \qquad \quad dV_x-V_{xx}dx-e^{V+W}dy,    dW_y-e^{V+W}dx-W_{yy}dy \rangle_{\text diff}   .
\end{aligned}
\EqTag{ECD5}
\end{equation}

For example the reduction $\CalI/\Gamma_{G_2}$ can be checked by noting
\begin{equation}
\begin{aligned}
&&\bfq^*_{\Gamma_{G_2}} (du-u_xdx-u_y dy) &= \  \frac{-2}{v+w} (\theta^v+\theta^w) +\frac{1}{v_x}\theta^v_x+\frac{1}{w_y} \theta^w_y,\\
&&\bfq^*_{\Gamma_{G_2}} (du_x-u_{xx}dx-e^u dy) &=\  \frac{2v_x}{(v+w)^2} (\theta^v+\theta^w) -\frac{2v_x^2+v_{xx}(v+w)}{v_x^2(v+w)}\theta^v_x+\frac{1}{v_x} \theta^v_{xx},\\
&&\bfq^*_{\Gamma_{G_2}} (du_y-e^udx-u_{yy} dy) &=\  \frac{2w_y}{(v+w)^2} (\theta^v+\theta^w) -\frac{2w_y^2+w_{yy}(v+w)}{w_y^2(v+w)^2}\theta^w_{y}+\frac{1}{w_y} \theta^w_{yy} ,
\end{aligned}
\end{equation}
or by using a cross-section $\sigma:M/G_2\to M$ (Lemma \StRef{PB3}). We also note that
\begin{equation}
\begin{aligned}
&&\bfp^*_2 (du-u_xdx-u_y dy) &=\ \theta^V + \theta^W  ,   \\
&&\bfp^*_2(du_x-u_{xx}dx-e^u dy) &=\ \theta^V_x-e^V \theta^V , \\
&&\bfp^*_2 (du_y-e^udx-u_{yy} dy) &= \ \theta^W_y -e^W \theta^W ,
\end{aligned}
\end{equation}
so that
$$
\CalI/\Gamma_H = \bfp_2^*( \CalI/\Gamma_{G_2} )+ \spn \{ \theta^ V \} \  .
$$
The structure equation 
$$
d \theta^V =- \theta^V_x \wedge dx + e^W  \theta^W \wedge dy,
$$
demonstrates Theorem \StRef{ORD}  by directly verifying $\CalI/\Gamma_H$ is an integral extension of $\CalI/\Gamma_{G_2}$. 

Finally we point out that the reductions in \EqRef{ECD5} are the standard Pfaffian systems for the partial differential equations
$$
z_{xy}=0, \quad u_{xy}=e^u, 
$$
and the prolongation of 
\begin{equation}
V_{y} =- e^W, W_x = -e^V . 
\EqTag{BTVW}
\end{equation}
Note that eliminating $V,W,V_x,W_y$ using equations \EqRef{DP12} gives
$$
z_x-u_x=  \sqrt{2} e^{\frac{z+u}{2}} \quad z_y+w_y = -\sqrt{2}e^{\frac{u-z}{2}} 
$$
which is the classical B\"acklund transformation relating the wave equation and Liouville's equation $u_{xy}=e^u$.
\end{Example}

\section{Iterated Reduction and De-prolongation}

\subsection{Iterated Reduction}

Let $\mu:G\times M \to M$ be a regular action of the Lie group $G$ on $M$ which preserves the regular integrable distribution $\Delta$,
\begin{equation}
(\mu_g)_* \Delta_x = \Delta_{\mu(g,x)} \, \quad {\rm for \ all \ } \ x \in M, \ g \in G.
\EqTag{GsymD}
\end{equation}
These assumption of regularity on the action of $G$ and the distribution $\Delta$ give rise to the two quotients
$$
\begin{gathered}
\beginDC{\commdiag}[3]
\Obj(-25, 0)[I]{$M$}
\Obj(0, 0)[H]{$M /G$}
\Obj(-25, -17)[I1]{$M /\Delta$ }
\mor{I}{H}{$\bfq_G$}[\atleft, \solidarrow]
\mor{I}{I1}{$\bfq_{\Delta}$}[\atright, \solidarrow]
\enddc
\end{gathered}
$$
where the functions $\bfq_G$ and $\bfq_{\Delta}$ are surjective submersions.

We now show, by virtue of the symmetry condition \EqRef{GsymD}, that the group $G$ acts naturally on $M/\Delta$ in which case we can construct the first iterated reduction $(M/\Delta)/G$. We begin with the following lemma.

\begin{Lemma} \StTag{gSp} Let $x\in M$, $g\in G$ and let $S_x$ be the maximal integral manifold of $\Delta$ through $x$. Then
$\mu(g, S_x) = S_{\mu(g,x)}$ where $S_{\mu(g,x)}$ is the maximal integral manifold of $\Delta$ through $\mu(g,x)$.
\end{Lemma}
\begin{proof}  Condition \EqRef{GsymD} implies, for all $g\in G$, that the function $\mu_g:M \to M$ maps integral manifolds of
$\Delta$ to integral manifolds. The image $\mu(g, S_x)$ is a connected integral manifold through the point $\mu(g,x)$ and is therefore contained in $S_{\mu(g,x)}$. Similarly $\mu(g^{-1}, S_{\mu(g,x)}) \subset S_x $. These two inclusions give
$$
S_{\mu(g,x)} = \mu(g ,\mu(g^{-1}, S_{\mu(g,x)} )) \subset \mu(g ,S_x)) \subset S_{\mu(g,x)} \,
$$
and so  $\mu(g, S_x) = S_{\mu(g,x)}$ as required.
\end{proof}

\begin{Corollary} \StTag{Cred} The function $\tilde \mu: G \times M/\Delta \to M/\Delta$
\begin{equation}
\tilde\mu(g, S_x)= S_{\mu(g,x)},  \quad g \in G, \ x \in M,
\EqTag{GC}
\end{equation}
is well-defined and defines an action of $G$ on $M/\Delta$. The projection map $\bfq_\Delta:M \to M/\Delta$ is equivariant with respect to the actions $\mu$ and $\tilde \mu$ of $G$, that is,
\begin{equation}
\bfq_\Delta ( \mu(g,x))  =\tilde \mu(g,  \bfq_\Delta(x)) , \quad {\rm for \ all \ } g \in G, \ x \in M .
\EqTag{GC2}
\end{equation}
\end{Corollary}

Let $\tilde \bfq_G : M/\Delta \to (M/\Delta)/G$ be the quotient map for the action in \EqRef{GC}.

\medskip

To define the second iterated reduction, we begin by using the $G$-invariance of $\Delta$ in equation \EqRef{GsymD} to define the bundle reduction  $\barDelta =\Delta/G\subset T(M/G)$ by
\begin{equation}
\barDelta = \Delta/G = \bfq_{G*}(\Delta) .
\EqTag{CQ}
\end{equation}
If the distribution $\barDelta$ on $M/G$ is integrable, then the second iterated reduction is $\bfq_{\barDelta}:M/G \to (M/G)/\barDelta $.

Theorem \StRef{Itred} given below states that if the action of $G$ on $M/\Delta$ is regular, then $\barDelta=\Delta/G$ is integrable and regular, and that the two iterated reductions $(M/\Delta)/G$ and $(M/G)/\barDelta$ are canonically diffeomorphic.  By an abuse of notation we then denote $\bfq_{\bar \Delta}:M/G \to (M/\Delta)/G$ which gives rise to the commutative diagram,
\begin{equation}
\begin{gathered}
\beginDC{\commdiag}[3]
\Obj(-25, 0)[I1]{$M$}
\Obj(0, 0)[I1G]{$M/G$}
\Obj(-25, -17)[I1C]{$M/\Delta$}
\Obj(0, -17)[I1Z]{$(M/\Delta)/G \, .$}
\mor{I1}{I1G}{$\bfq_G$}[\atleft, \solidarrow]
\mor{I1G}{I1Z}{$\bfq_{\barDelta}$}[\atleft, \solidarrow]
\mor{I1C}{I1Z}{$\tilde\bfq_G$}[\atright, \solidarrow]
\mor{I1}{I1C}{$\bfq_\Delta$}[\atright, \solidarrow]
\mor{I1}{I1Z}{$\pi$}[\atright, \solidarrow]
\enddc
\end{gathered}
\EqTag{CGreduce}
\end{equation}
This is the {\deffont iterated reduction diagram} for manifolds with $G$ invariant integrable distributions.

We now present two lemmas which are essential for the proof of the Theorem \StRef{Itred} (iterated reduction).

\begin{Lemma} \StTag{pmap}  There exists a unique function  $\bfp:M/G \to (M/\Delta)/G$, such that the diagram
\begin{equation}
\begin{gathered}
\beginDC{\commdiag}[3]
\Obj(-25, 0)[I]{$M$}
\Obj(0, 0)[I1]{$M/G$}
\Obj(0, -17)[I2]{$(M/\Delta)/G$}
\Obj(-25, -17)[I3]{$M/\Delta$}
\mor{I}{I1}{$\bfq_G$}[\atleft, \solidarrow]
\mor{I1}{I2}{$\bfp$}[\atleft, \solidarrow]
\mor{I}{I3}{$\bfq_{\Delta}$}[\atright, \solidarrow]
\mor{I3}{I2}{$\tilde \bfq_G$}[\atright, \solidarrow]
\mor{I}{I2}{$\pi$}[\atleft, \solidarrow]
\enddc
\EqTag{defp}
\end{gathered}
\end{equation}
commutes, where $\pi = \tilde \bfq_G \circ \bfq_\Delta$.
\end{Lemma}

\begin{proof} Let $\mu(G,x) \in M/G$, and let
\begin{equation}
\bfp(\mu(G,x))  =  \tilde \bfq_G \circ \bfq_{\Delta}(x).
\EqTag{bfpd}
\end{equation}
We need to check that $\bfp$ in \EqRef{bfpd} is well-defined by showing by showing $\tilde \bfq_G \circ \bfq_{\Delta}(x)$ is
independent of the point $x$ in the orbit $\mu(G,x)$.  Let $x,x' \in M$ be two points which satisfy  $\bfq_G(x) = \bfq_G(x')$. If $\tilde \bfq_G \circ \bfq_\Delta(x)= \tilde \bfq_G \circ \bfq_\Delta(x')$ then the function $\bfp$ will be well defined as a function on $M/G$.

Since $\bfq_G(x)= \bfq_G(x')$, then $x'= \mu(g,x)$ for some $g\in G$. We then compute,
$$
\begin{aligned}
\tilde \bfq_G \circ \bfq_\Delta(x') & =\tilde \bfq_G(\Delta(\mu(g,x)) \\
& = \tilde \bfq_G \circ \tilde \mu(g, \bfq_\Delta(x)), \qquad {\rm by \ equivariance \ in \ equation \ \EqRef{GC2}} \\
& = \tilde \bfq_G\circ \bfq_\Delta (x).
\end{aligned}
$$
Therefore $\bfp$ in equation \EqRef{bfpd} is well-defined.
\end{proof}

We now give the second lemma.

\begin{Lemma} \StTag{RmodG} Let $\bfp:N \to Q$ be a surjective submersion with connected fibres. Then $\ker (\bfp_*)$ is
 an integrable and regular distribution, and $N/\ker(\bfp_*)$ is canonically diffeomorphic to $Q$.
\end{Lemma}
\begin{proof} The distribution $\ker (\bfp_*)$ is integrable and so we need to show it is regular.

We begin with a standard construction (see for example \cite{munkres:2000a}, page 142). Define a partition on $N$ by $x \mytilde x'$ if $\bfp(x) = \bfp(x')$. The equivalence class of $x\in N$ is then
\begin{equation}
[x] = \bfp^{-1} (\bfp(x))\, .
\EqTag{pretilde}
\end{equation}
The function $\bfp$ is a submersion and therefore an open map, and the canonical bijection $
\tau :Q  \to N/\!\!\mytilde $ given by
\begin{equation}
\tau(y) = \bfp^{-1}(y)
\EqTag{pretau}
\end{equation}
is a homeomorphism where $N/\!\!\mytilde$ has the quotient topology. The differentiable structure on $Q$ induces a differentiable structure on $N/\!\!\mytilde$ compatible with the quotient topology in which case the map $\tau$ is a diffeomorphism. With this differentiable structure on $N/\!\!\mytilde$,  the function $\tilde \bfq:N \to N/\!\!\mytilde$ is then a surjective submersion.

We now show that $\ker (\bfp_*)$ is regular by first showing that for each $x \in N$,  $[x]= \bfp^{-1}(\bfp(x)) = S_x$, where $S_x$
is the maximal connected integral manifold of $\ker (\bfp_*)$ through $x$. Therefore the elements of $N/\!\!\mytilde$ and $N/\ker( \bfp_*)$ are
identical, the functions $\bfq_{\ker (\bfp_*)} $, and $ \tilde \bfq$ are identical and so $\bfq_{(\ker \bfp_*)}$ is a surjective submersion.

Let $x\in N$. Then $\bfp^{-1}(\bfp(x))$ is a connected (by hypothesis) integral manifold of $\ker( \bfp_*)$ and therefore $[x] = \bfp^{-1}(\bfp(x)) \subset S_x$. To show the reverse inclusion let $x' \in S_x$,
and since $S_x$ is connected, let $\gamma:[0,1] \to S_x$ be a smooth curve in $S_x$ joining $x$ to $x'$. Since $\gamma(t) \in S_x$, and $\gamma(t)$ is smooth as a curve into $N$ we have the tangent vector $\dot \gamma(t)$ satisfying $\dot \gamma(t) \in \ker (\bfp_*)$
for all $t\in (0,1)$. Consequently, the tangent vector to the curve $\bfp \circ \gamma: [0,1] \to Q$  satisfies $ \frac{d\ }{dt}(\bfq \circ \gamma)= 0$, for all $t\in (0,1)$. Therefore $\bfp \circ \gamma$ is constant and $\bfp(x')=\bfp\circ\gamma(1)=\bfp\circ\gamma(0)= \bfp(x)$. That is $S_x \subset \bfp^{-1}(\bfp(x))$, which then shows $S_x = \bfp^{-1}(\bfp(x))$. Therefore the elements of $N/\!\!\mytilde$ and $N/\ker( \bfp_*)$ are identical which finishes the proof.
\end{proof}

We now come to the main theorem on iterated reduction.

\begin{Theorem} \StTag{Itred} Let $G$ be a Lie group preserving the regular integrable distribution $\Delta$ on the manifold $M$, and suppose that $G$ acts regularly on $M$ and  on $M/\Delta$. Then
\noindent
{\bf [i]} the distribution $\barDelta=\Delta/G $ on $M/G$ is integrable and regular, and
\noindent
{\bf [ii]} there exists a canonical diffeomorphism of $\tau:(M/\Delta)/G \to (M/G)/\barDelta$, giving rise to the
commutative diagram
\begin{equation}
\begin{gathered}
\beginDC{\commdiag}[3]
\Obj(-25, 0)[I]{$M/G$}
\Obj(0, -17)[I2]{$(M/G)/{\bar \Delta}$}
\Obj(-25, -17)[I3]{$(M/\Delta)/G$}
\mor{I}{I2}{$\bfq_{\bar \Delta}$}[\atleft, \solidarrow]
\mor{I}{I3}{$\bfp$}[\atright, \solidarrow]
\mor{I2}{I3}{$\tau^{-1}$}[\atleft, \solidarrow]
\enddc
\EqTag{mytauinv}
\end{gathered}
\end{equation}

\end{Theorem}

\begin{proof} The hypothesis that $G$ acts regularly on $M/\Delta$ implies that $\pi =\tilde \bfq_G \circ \bfq_{\Delta}$ in diagram \EqRef{defp} is a surjective submersion. By applying a generalization of Theorem \StRef{ORD} (to $\pi$ and $\bfq_G$ in diagram \EqRef{defp}) we have $\bfp$ is a surjective submersion (see Theorem 3.1 in \cite{anderson-fels:2011a}).

We now show $\ker (\bfp_*)$ is integrable and regular and that $\ker(\bfp_*)=\barDelta$. Let $\tilde x \in (M/\Delta)/G$, then $\bfp^{-1}(\tilde x) = \bfq_G( S_x)$ where $x\in M$ and $S_x$ is the maximal integral manifold of $\Delta $ through $x\in M$. In particular $\bfp^{-1}(\tilde x)$ is connected. By Lemma \StRef{RmodG}, $\ker (\bfp_*)$ is integrable and regular. Now from the commutative diagram \EqRef{defp} we have
$$
\ker (\bfp_*)  = \bfq_{G*}( \ker( \pi_*) )  = \bfq_{G*}(\Delta+\bfGamma_G) = \bfq_{G*} (\Delta) = \Delta/G  .
$$
Therefore  $\barDelta=\Delta/G$ is an  integrable regular distribution. This proves part {\bf [i]}.

We now apply Lemma \StRef{RmodG} with $N=M/G$, $Q=(M/\Delta)/G$ and $\ker (\bfp_*) = \barDelta$. The canonical diffeomorphism $\tau:(M/\Delta)/G \to (M/G)/\barDelta$ in equation \EqRef{pretau} is then
\begin{equation}
\tau(\mu(G, S_x)) = \bar S_{\mu(G,x)}\, ,
\EqTag{mytau}
\end{equation}
where $\bar S_{\mu(G,x)}$ is the maximal integral manifold of $\barDelta$ through $\mu(G,x) \in M/G$. This proves part {\bf [ii]} of the theorem.
\end{proof}

Theorem \StRef{Itred} provides the identification of $(M/\Delta)/G$ and $(M/G)/\barDelta$. In terms of this identification the function $\bfp $ in equation \EqRef{defp} becomes $\bfq_{\barDelta}$ so that diagram \EqRef{defp} becomes diagram \EqRef{CGreduce}.

\begin{Remark}  Let $\bfp_1:M \to (M/\Delta)/G$ be the surjective submersion $\bfp_1= \tilde \bfq_G \circ \bfq_\Delta$. As described in the proof of Lemma \StRef{RmodG}, the manifold $(M/\Delta)/G$ is canonically diffeomorphic to the quotient space $M/\!\!\mytilde_1$, where the $\mytilde_1$ equivalence class $[x]_1$ of $x\in M$ is given by
\begin{equation}
[x]_1 = \bfp_1^{-1}( \bfp_1(x)) = \mu(G, S_x)\,
\EqTag{myt1}
\end{equation}
and $S_x$ is the maximal integral manifold of $\Delta$ through $x$.

Similarly, let  $\bfp_2 :M \to (M/G)/\barDelta$ be the surjective submersion $\bfp_2(x) = \bfq_{\barDelta} \circ \bfq_G(x)$. Then  $(M/G)/\barDelta $ is canonically diffeomorphic to the quotient space $ M/\!\!\mytilde_2$,  where the $\mytilde_2$ equivalence class, $[x]_2$ of $x\in M$ is given by
\begin{equation}
[x]_2 =\bfp_2^{-1} (\bfp_2(x)) = \bfq_G^{-1}( \bar S_{\bfq_G(x)})\,
\EqTag{myt2}
\end{equation}
and $\bar S_{\bfq_G(x)} $ is the maximal integral manifold of $\barDelta$ through $\bfq_G(x)$.

Theorem \StRef{Itred} implies that the equivalence classes \EqRef{myt1} and \EqRef{myt2} are identical. By applying the projection $\bfq_G$ to these equivalence classes we find
\begin{equation}
\bfq_G([x]_1) = \bfq_G(\mu(G,S_x))=\bfq_G(S_x),  \quad {\rm for \ all \ }   x \in M.
\EqTag{ece1}
\end{equation}
and
\begin{equation}
\bfq_G([x]_2) =   \bar S_{\bfq_G(x)} , \quad {\rm for \ all \ }   x \in M.
\EqTag{ece2}
\end{equation}
In other words, by equations \EqRef{ece1} and \EqRef{ece2}, the equivalence classes \EqRef{myt1} and  \EqRef{myt2} are the same if and only if $\bfq_G$ maps maximal integral manifolds of $\Delta$ to maximal integral manifolds of the quotient $\barDelta$,
\begin{equation}
\bfq_G(S_x) =  \bar S_{\bfq_G(x)}  \quad {\rm for \ all \ }   x \in M.
\EqTag{GSx}
\end{equation}
Condition \EqRef{GSx} is precisely the condition we need in order to construct the commutative diagram \EqRef{CGreduce} and is satisfied given the regularity hypothesis stated in Theorem \StRef{Itred}.
\end{Remark}

\subsection{De-prolongation}

Recall from section \StRef{S23} that  the Cauchy characteristics $\Cau(\CalI) \subset TM$ of an EDS  $\CalI$ is an integrable distribution which we will assume to be regular where $\bfq_{\Cau_{\CalI}}: M \to M/\Cau_{\CalI} $ is the quotient map and  $\bar \CalI = \CalI/\bfq_{\Cau_{\CalI}}$  the reduced EDS (see Theorem \StRef{CQ}).  Combining iterated reduction and Cauchy reduction produces the following theorem.

\begin{Theorem}\StTag{SymmodCI} Let $G$ be a symmetry group of $\CalI$ acting regularly on $M$. Then
\\
{\bf [i]} the Cauchy characteristics $\Cau(\CalI) $ are  $G$ invariant and,
\\
{\bf [ii]} $G$ is a symmetry group of  $\CalI/\bfq_{\Cau_{\CalI}} $, where the action of $G$ on $M/\Cau_{\CalI}$ is defined by equation \EqRef{GC} (with $\Delta = \Cau(\CalI)$).
\\
{\bf [iii]} If the action of $G$ on $M$ and $M/\Cau_{\CalI}$ is regular then the following diagram commutes,
\begin{equation}
\begin{gathered}
\beginDC{\commdiag}[3]
\Obj(0, 0)[I1]{$\left(\, M,\ \CalI\, \right)$}
\Obj(50, 0)[I1G]{$\left(\, M/G,\ \CalI/G\, \right)$}
\Obj(0, -15)[I1C]{$\left(\, M/{\Cau_{\CalI}},\ \CalI/\bfq_{\Cau_{\CalI}}\, \right)$}
\Obj(50, -15)[I1Z]{$\left(\, (M/{\Cau_{\CalI}})/G,\ (\CalI/\bfq_{\Cau_{\CalI}})/G\, \right)\ ,$}
\mor{I1}{I1G}{$\bfq_G$}[\atleft, \solidarrow]
\mor{I1G}{I1Z}{$\bfq_{{\bar \Cau_{\CalI}}}$}[\atleft, \solidarrow]
\mor{I1C}{I1Z}{$\tilde\bfq_G$}[\atright, \solidarrow]
\mor{I1}{I1C}{$\bfq_{\Cau_{\CalI}}$}[\atright, \solidarrow]
\mor(3,-2)(37,-13){$\pi$}[\atright, \solidarrow]
\enddc
\end{gathered}
\EqTag{CGreduce2}
\end{equation}
where $\tilde \bfq_G$ is the quotient by the action of $G$ on $M/\Cau_{\CalI}$ from part  {\bf [ii]}, and $\bar \Cau_{\CalI}= \Cau(\CalI)/G$.
\end{Theorem}

\begin{proof} We begin by proving $\Cau_{\CalI}$ is $G$ invariant. If $X \in \Cau(\CalI)$ then
$$
g_* X \hook \CalI = X \hook (g^* \CalI) = X \hook \CalI \subset \CalI .
$$
Therefore $g_*X$ is a Cauchy characteristic and $\Cau(\CalI)$ is $G$-invariant. This proves part {\bf [i]}.

We now check that the action of $G$ on $M/\Cau_{\CalI}$ is a symmetry of $\CalI/\bfq_{\Cau_{\CalI}}$. Let $ \bar \theta \in \CalI/\bfq_{\Cau_{\CalI}}$, and $g\in G$. We compute $\bfq_{\Cau_{\CalI}}^*\circ \tilde \mu_g^* (\bar \theta)$ using the equivariance of $\bfq_{\Cau_{\CalI}}$ (equation \EqRef{GC2}) to obtain
\begin{equation}
\bfq_{\Cau_{\CalI}}^*\circ \tilde \mu_g^* (\bar \theta) = \mu_g^*\circ \bfq_{\Cau_{\CalI}}^* (\bar \theta) .
\EqTag{GonI}
\end{equation}
By hypothesis $\bfq_{\Cau_{\CalI}}^* (\bar \theta) \in \CalI$ and $G$ is a symmetry group of $\CalI$, therefore $ \mu_g^* \circ \bfq_{\Cau_{\CalI}}^* (\bar \theta) \in \CalI$. Together with equation \EqRef{GonI}, this shows $\tilde \mu_g^* (\bar \theta) \in \CalI/\bfq_{\Cau_{\CalI}}$, and thus $G$ is a symmetry group of $\CalI/\bfq_{\Cau_{\CalI}}$. This proves part {\bf [ii]}.

We now prove part {\bf [iii]}. Part {\bf [i]} shows that $\Cau(\CalI)$ is $G$ invariant and, if we assume that the action of $G$ is regular on both $M$ and $M/\Cau(\CalI)$, then Theorem \StRef{Itred} can be applied to construct the commutative diagram \EqRef{CGreduce}, again with $\Delta = \Cau(\CalI)$ and $\bar \Cau_{\CalI} = \Cau(\CalI)/G$. Combining this with a simple generalization of Theorem \StRef{ORD} (see Theorem 3.1 in \cite{anderson-fels:2011a}) we obtain the commutative diagram of differential systems in \EqRef{CGreduce2}.
\end{proof}

The issue we now address is to what extent the reduced distribution $\bar\Cau_{\CalI}=\Cau(\CalI)/G$ in \EqRef{CGreduce2} are precisely the Cauchy characteristics of $\CalI/G$. Denote the Cauchy characteristics of the reduced system $\CalI/G$ by $\Cau_{\CalI/G}=\Cau(\CalI/G)$.  If diagram \EqRef{CGreduce} holds and if the condition $\Cau(\CalI) /G = \Cau(\CalI/G)$ is satisfied, then we say that {\it characteristic reduction and symmetry reduction commute}. We investigate this condition starting with the following lemma.

\begin{Lemma} \StTag{CGlemma} Let $\Cau_{\CalI} \cap \bfGamma_G$ be constant rank. Then the distribution $\bar \Cau_{\CalI} =\Cau_{\CalI}/G= \bfq_{G*}( \Cau_{\CalI})$ has constant rank, is integrable, and  ${\Cau_{\CalI}}/G \subset \Cau(\CalI/G)$.
\end{Lemma}
\begin{proof} The proof that $\bar \Cau_{\CalI}$ is constant rank and integrable follows from Section 2.5 and is also given by Theorem 2.9 in \cite{fels:2007a}.  We proceed with the second part of the lemma.

Let $ y \in M/G$, $Y \in (\bar \Cau_{\CalI})_y$, $\bar \theta \in (\CalI/G)_{y}$, and $x\in M$ with $\bfq_G(x) = y$. To show $Y \in  \Cau(\CalI/G)_y$, we need to show $Y \hook \bar \theta \in (\CalI/G)_y$ which is equivalent to showing $\bfq_G^*( Y\hook \bar \theta) \in \CalI_x$. Now since $Y\in (\bar \Cau_{\CalI})_y$ we may choose $X \in (\Cau_{\CalI})_x$ satisfying $\bfq_{G*}( X ) = Y$. We then have
\begin{equation}
\bfq_G^*( Y \hook \bar \theta) = X \hook \bfq_G^* \bar \theta \, .
\EqTag{CC1}
\end{equation}
Now $X$ is a characteristic of $\CalI$ and so $X \hook \bfq_G^* \bar \theta \in \CalI_x$. Therefore by equation \EqRef{CC1}, $Y\hook \bar \theta \in (\CalI/G)_y$, which proves the lemma.
\end{proof}

Lemma \StRef{CGlemma} applied to diagram \EqRef{CGreduce2} allows us to conclude that symmetry reduction and Cauchy reduction commute with the following caveat - {\it reduction by $ \bar \Cau_{\CalI} = \Cau(\CalI)/G$ in diagram \EqRef{CGreduce2} may only be by a subset of the Cauchy characteristics of $\CalI/G$}. In the case where the action of $G$ is transverse to $\CalI$ (see equation \EqRef{Itrans}) more can be said.

\begin{Theorem}\StTag{Cred2} Let $G$ be a symmetry group of the EDS $\CalI$ with Cauchy characteristics $\Cau_{\CalI}$. Suppose $\bfGamma _G\cap \Cau_{\CalI}$ is constant rank and that $G$ acts transversally to $\CalI$. Then the reduced characteristics  $\bar \Cau_{\CalI}=\Cau_{\CalI}/G$ are the Cauchy characteristics for the reduced EDS $\CalI/G$, that is
$$
\Cau(\CalI)/G= \Cau(\CalI/G), \quad {\text or} \quad (\Cau_{\CalI})/G= \Cau_{\CalI/G} .
$$
\end{Theorem}
\begin{proof} By Lemma \StRef{CGlemma} we have $\Cau(\CalI)/G=\bfq_{G*}(\Cau_{\CalI}) \subset \Cau(\CalI/G)$, and so we just need to show the reverse inclusion. We utilize the transversality hypothesis by working with a set of local generators for $\CalI$ and $\CalI/G$  chosen such that (see \cite{anderson-fels:2011a},  \cite{anderson-fels:2005a}),
\begin{equation}
	\CalI|_{U}  = \langle\, \theta^a,\,  \theta^i_{G}, \,  \tau^\alpha_{G} \, \rangle_\text{alg}, \quad {\rm and} \quad
(\CalI/G)|_{\bfq_G(U)}  = \langle\,   \bar \theta^i, \,  \bar \tau^\alpha \, \rangle_\text{alg},
\EqTag{CC10}
\end{equation}
where
\begin{equation}
\bfq_G^* ( \bar \theta^i) = \theta^i_G, \quad \bfq_G^* (\bar \tau^\alpha) = \tau^\alpha_G.
\EqTag{CC11}
\end{equation}

Let $y \in \bfq_G(U)$ and suppose $Y \in \Cau(\CalI/G)_{y}$ is a Cauchy characteristic at $y$. Pick $x\in U$ with $\bfq_G(x) = y$ and $X\in T_xU$ with $ \bfq_{G*}( X) = Y$. The transverallity hypothesis implies there exists a local basis $\{Z_b\}$ for $ \bfGamma_G$ at $x$ such that
\begin{equation}
\theta^a(Z_b) = \delta^a_b .
\EqTag{defZa}
\end{equation}
We claim that the tangent vector,
\begin{equation}
Z= X -\theta^a(X)Z_a,
\EqTag{myZ}
\end{equation}
which satisfies $(\bfq_G)_* Z = Y$, also satisfies $Z \in \Cau(\CalI)_x$. Once this is shown, the theorem is proved.

We check the characteristic condition for $Z$ on the algebraic generators of $\CalI|_U$ in equation \EqRef{CC10}. Clearly, by equations \EqRef{defZa} and \EqRef{myZ}, $ \theta^a(Z) = 0 $. Now if $\omega_G$ is a $G$-basic form in $\CalI$, then $\omega_G=\bfq_G^* \bar \omega$, where $\bar \omega \in \CalI/G$ and
\begin{equation}
\begin{aligned}
Z \hook \omega_G & = X \hook \bfq_G^* \bar \omega & & \ {\rm using }  \omega_G \ {\rm is \ semi-basic } &\\
& =\bfq_G^* \left((\bfq_G)_* X \hook \bar \omega \right) & &  & \\
& =\bfq_G^* \left(Y \hook \bar \omega \right),& &\ {\rm using \ } \bfq_{G*}(X) = Y \quad . \ &
\end{aligned}
\EqTag{CC20}
\end{equation}
By hypothesis, $Y$ is a Cauchy characteristic and so $Y \hook \bar \omega \in (\CalI/G)_y$. The last line in equation \EqRef{CC20} then implies $Z \hook \omega_G \in \CalI$. Applying this result to the generators $\theta^i_G$ and $\tau^\alpha_G$ in $\CalI|_U$ in equation \EqRef{CC10} implies
 that $Z \in \Cau(\CalI)_x$ and is a Cauchy characteristic of $\CalI$.
\end{proof}

\begin{Remark} \StTag{COMRD} Theorem \StRef{Cred2} shows for an action of $G$ transverse to $\CalI$, that symmetry reduction and reduction by Cauchy characteristics commute. In which case we may replace $\bfq_{\bar \Cau_{\CalI}}$ by $\bfq_{\Cau_{\CalI/G}}$ in diagram \EqRef{CGreduce2} since $\Cau(\CalI/G)= \Cau(\CalI)/G$. 
\end{Remark}

Finally, we apply iterated reduction to the situation which will allow us to pass from diagram \EqRef{Intro10} to diagram \EqRef{FCC}. Let $\CalI$ be a Pfaffian system on $M$, with symmetry group $G$ acting transversally to the derived system $\CalI'$. Then $\CalI/G$ is a Pfaffian system \cite{anderson-fels:2005a}. However, even with the  hypothesis that $G$ acts transversally to $\CalI'$, it is not necessarily true that $(\CalI/G)' = \CalI'/G$. One easy way to see why these may not be equal is
that $(\CalI/G)'$ is always a Pfaffian system, while $\CalI'/G$ is a Pfaffian system if and only if $G$ is transverse to $\CalI''$.

Let $ \Cau(\CalI') $ be the Cauchy characteristics of $\CalI'$.  The analogue to Theorem \StRef{SymmodCI} is the following.

\begin{Theorem}\StTag{SymmodC} Let $G$ be a symmetry group of the Pfaffian system $\CalI$. 
\\
{\bf [i]} The Cauchy characteristics $\Cau(\CalI')$ of $\CalI'$ are  $G$ invariant and,
\\
{\bf [ii]} $G$ is a symmetry group of  $\CalI/\bfq_{\Cau_{\CalI'}} $, where the action on $M/\Cau(\CalI')$ is defined by equation \EqRef{GC} (with $\Delta = \Cau(\CalI')$).
\\
{\bf [iii]} If $G$ acts regularly on $M$ and $M/\Cau_{\CalI'}$  then the following diagram commutes
\begin{equation}
\begin{gathered}
\beginDC{\commdiag}[3]
\Obj(0, 32)[I1]{$(M,\CalI)$}
\Obj(0, 18)[I1G]{$(M/\Cau_{\CalI'},\CalI/\bfq_{\Cau_{\CalI'}})$}
\Obj(50, 32)[I1C]{$(M/G,\CalI/G)$}
\Obj(50, 18)[I1Z]{$( (M/\Cau_{\CalI'})/G,(\CalI/ \bfq_{\Cau_{\CalI'}})/G)$}
\mor{I1}{I1G}{$\bfq_{\Cau_{\CalI'}}$}[\atleft, \solidarrow]
\mor{I1G}{I1Z}{$\tilde \bfq_G$}[\atleft, \solidarrow]
\mor{I1C}{I1Z}{$\bfq_{\bar \Cau_{\CalI'}}$}[\atright, \solidarrow]
\mor{I1}{I1C}{$\bfq_{G}$}[\atright, \solidarrow]
\mor(3,30)(42,20){$\pi$}[\atright, \solidarrow]
\enddc
\end{gathered}
\EqTag{Ddepro}
\end{equation}
where $\tilde \bfq_{\tilde G}$ is the quotient by the action of $G$ on $M/\Cau_{\CalI'}$.
\\
{\bf [iv]} If $G$ is transverse to $\CalI'$ then $\Cau(\CalI')/G = \Cau(\CalI'/G) $ and
\\
{\bf [v]} $ \Cau(\CalI'/G) \subset \Cau\left( (\CalI/G)' \right) $, where $\Cau \left( ( \CalI/G)' \right)$ are the Cauchy characteristics of the Pfaffian system $(\CalI/G)'$.
\end{Theorem}

\begin{proof}
The proofs of {\bf [i],[ii],[iii]}  are the same as in Theorem \StRef{SymmodCI}.  Part {\bf [iv]} is a direct consequence of Theorem \StRef{Cred2} applied to $\CalI'$ and so there is nothing to prove.

To prove part {\bf [v]} let $Y\in \Cau(\CalI'/G) $. We show that $Y \in \Cau \left( ( \CalI/G)' \right)$ by
checking the Cauchy characteristic condition on generators for $(\CalI/G)'$. We begin with the 1-forms in $\CalI'/G$ and $(\CalI/G)'$ which by equation \EqRef{Qdone} in Corollary \StRef{RDone} are the same. Therefore the Cauchy characteristic condition on the one-forms for these two systems is identical. Since $(\CalI/G)'$ is Pfaffian we only need to check that if $\bar \tau \in (\CalI/G)' \cap\, \Omega^2(M/G)$
then $ Y \hook \bar \tau  \in (\CalI/G)'$.

By Lemma \StRef{RDone} $(\CalI/G)' \subset \CalI'/G$, and so $\bar \tau \in (\CalI'/G)$. Since $Y$ is a characteristic for $\CalI'/G$,
$$
Y \hook \bar \tau  \in (\CalI'/G) \cap \Omega^1(M/G) .
$$
By equation \EqRef{Qdone} $(\CalI'/G) \cap \Omega^1(M/G) = (\CalI/G)' \cap \Omega^1(M/G)$, and so $Y\hook \bar \tau \in (\CalI/G)'$. Therefore $Y \in \Cau(\, (\CalI/G)' \, )$. 
\end{proof}

\begin{Remark}\StTag{CPD} In Theorem \StRef{SymmodC} the set  $ \Cau({\CalI'})/G$ is not necessarily the full space of Cauchy characteristics $\Cau((\CalI/G)')$. However, suppose that the conditions in Theorem  \StRef{SymmodC} hold and 

{\bf [i]}\  $pr( \CalI/\bfq_{\Cau_{\CalI'}} ) = \CalI \quad $ and 

{\bf [iI]} \ $\Cau(\CalI')/G = \Cau( (\CalI/G)' )\quad $ 

\noindent
where $pr$ is prolongation. Then $\CalI/\bfq_{\Cau_{\CalI'}}$ is the {\deffont de-prolongation  of} $\CalI$,
the quotient by $\Cau(\CalI')/G$ in diagram \EqRef{Ddepro} can be replaced with the quotient by $\Cau\left( (\CalI/G)' \right)$, and the diagram \EqRef{Ddepro} is the {\deffont de-prolongation diagram}.  With these additional hypothesis, {\deffont de-prolongation and group reduction commute}. 
\end{Remark}

\begin{Example} As a precursor to the proofs of Theorems \StRef{BTMA0} and \StRef{BTMA}, and to demonstrate de-prolongation,  we revisit Example \StRef{EX1}.  In this case we have in equation \EqRef{ECD5} the rank 3 Pfaffian systems $\CalI/\Gamma_{G_a}$
on the 7 manifolds $M/\Gamma_{G_a}$  and the rank 4 Paffian system $\CalB=\CalI/\Gamma_H$ on the 8 manifold $M/\Gamma_H$. 
Now from equation \EqRef{Twof} we have
$$
\Cau(\CalI') = \spn \{\ \partial_{v_{xxx}}, \ \partial_{w_{yyy}} \ \} \, 
$$
and $M/\Cau(\CalI') =J^2(\real,\real)\times J^2(\real,\real)$, and $\CalI/\Cau(\CalI') = \CalC_1+\CalC_2$ (the standard contact structures). The reductions $\Cau(\CalI')/\Gamma_{G_a}$ and $\Cau(\CalI')/\Gamma_{H}$ are all two dimensional and satisfy
\begin{equation}
\begin{aligned}
\Cau(\CalI')/\Gamma_{G_1}\, &=\, \spn \{ \partial_{z_{xx}}, \partial_{z_{yy}} \}  &\!\!\!\!&  =\, \Cau_{\CalI_{s=0}'} \\
\Cau(\CalI')/\Gamma_{G_2}\, &=\, \spn \{ \partial_{u_{xx}}, \partial_{u_{yy}} \}  &\!\!\!\!&  =\, \Cau_{\CalI_{2}'}\quad  \\
\Cau(\CalI')/\Gamma_{H}\, &=\, \spn \{ \partial_{V_{xx}}, \partial_{W_{yy}} \}  &\!\!\!\!& =\,  \Cau_{\CalB'}\quad
\end{aligned}
\EqTag{ECDd}
\end{equation}
Therefore conditions {\bf [i],[ii]} in Remark \StRef{CPD} hold and we can form the de-prolongation diagram \EqRef{Ddepro} for each of $(\CalI,\Gamma_{G_1})$,  $(\CalI,\Gamma_{G_2})$, and $(\CalI, {\Gamma_H})$.

Let  $\tilde \Gamma_{G_a}$ and $\tilde \Gamma_H$ be the projection (or reduction by $\Cau(\CalI')$ as in Corollary \StRef{Cred}) of the Lie algebras $\Gamma_{G_a}$ and $\Gamma_H$ to $M/\Cau(\CalI')=J^2(\real,\real)\times J^2(\real,\real)$. 
Combing the application of Theorem \StRef{SymmodC} three times, we get the commutative diagram
\begin{equation*}
\beginDC{\commdiag}[3]
\Obj(0, 52)[I]{$ \boxed{\begin{gathered} dv - v_xdx, \ dv_x- v_{xx}\,dx,   \\[-2\jot] dw -w_y dy, \ dw_y-w_{yy}dy   \end{gathered}  \quad (\CalC_1+\CalC_2) } $ }
\Obj(0, 26)[H]{$\boxed{ \begin{gathered}  dV-V_x dx+e^W dy  ,\  \\[-2\jot] dW+e^Vdx-W_ydy        \end{gathered}  \quad  
(\CalB/\Cau_{\CalB' } )} $}
\Obj(-37, -1)[I1]{$\boxed{ \begin{gathered} dz - z_x dx-z_y dy \\[-2\jot]  dz_x \wedge dx,\ dz_{yy} \wedge dy  \end{gathered} \quad  \ (  \CalI_{s=0}/\Cau_{\CalI_{s=0}'}  \, ) }  $}
\Obj(37, -1)[I2]{$ \boxed{ \begin{gathered} du - u_x dx-u_y dy  \\[-2\jot]    (du_{x}  -e^u dy)\wedge dx,\ (du_y-e^u dx) \wedge dy  \end{gathered}   \quad  (\CalI_2/\Cau_{\CalI_2'} )   }$}
\mor{I}{H}{$\bfq_{\tilde \Gamma_H}$}[\atleft, \solidarrow]
\mor(-4, 17)(-45, 6){$\tilde \bfp_1$}[\atleft, \solidarrow]
\mor(4, 17)(42, 6){$\tilde \bfp_2$}[\atright, \solidarrow]
\mor(-30, 48)(-52, 4){$\bfq_{\tilde \Gamma_{G_1}}$}[\atright, \solidarrow]
\mor(30, 48)(52, 5){$\bfq_{\tilde \Gamma_{G_2} \ .}$}[\atleft, \solidarrow]
\enddc
\EqTag{ECD3}
\end{equation*}
where $\bfq_{\tilde \Gamma_{G_1}}, \bfq_{\tilde \Gamma_{G_1}}, \bfq_{\tilde \Gamma_{H}}$ are each the bottom arrow in the  corresponding diagram \EqRef{Ddepro}.

Diagram \EqRef{ECD3} produces the B\"acklund transformation in equation \EqRef{FCC} where $\tilde \bfp_1$ and $\tilde \bfp_2$ are the
de-prolongation of the maps $\bfp_1$ and $\bfp_2$ in equation \EqRef{DP12},  and are given in coordinates by
\begin{equation*}
\begin{aligned}
&\tilde \bfp_1= (x=x,y=y, z= V-W,  z_x = V_x+e^V,z_y = -W_y-e^W)\\
&\tilde \bfp_2 = (x=x,y=y, u=\log (2)+V+W, u_x = V_x- e^V , u_y= W_y- e^W )
\end{aligned}
\end{equation*}

\end{Example}

\section{B\"acklund Transformations for Darboux Integrable PDE in the Plane}

\subsection{The Vessiot Algebra and Quotient Representation}

Let $F(x,y,u,u_x,u_y, u_{xx},u_{xy},u_{yy})=0$ be a hyperbolic PDE in the plane and let $\CalI$ be the standard rank 3 hyperbolic Pfaffian system determined by the restriction of the contact system on $J^2(\real^2,\real)$ \cite{gardner-kamran:1993a} to the level set
$ M = \{\ p \in J^2(\real^2,\real) \ | \ F(p) = 0 \ \} $. The following theorem is proved in \cite{gardner-kamran:1993a}.

\begin{Theorem}\label{hyperpde} Let $ \CalI$ be the standard rank 3 Pfaffian system for a hyperbolic PDE in the plane on the seven dimensional manifold $M$. Then about each point $x\in M$ there exists an open set $U$ and  a coframe $\{ \theta^i, \omega^a,\pi^a\}_{0\leq i\leq 2, 1\leq a,\leq 2} $ on $U$  where $\CalI|_U=\langle \theta^0,\theta^1,\theta^2 \rangle_{\text diff}$ and 
\begin{equation}
\begin{aligned}
d\theta ^0 & = \theta^1 \wedge \omega ^1 + \theta^2 \wedge \omega ^2 & & \mod \theta^0\\
d \theta^1 & = \omega^1 \wedge \pi^1 + \mu_1 \theta^2 \wedge \pi^2 & &\mod \theta^0, \theta^ 1 \\
d \theta^2 & = \omega^2 \wedge \pi^2 + \mu_2 \theta^1 \wedge \pi^1 & &\mod \theta^0, \theta^ 2  \  .
\end{aligned}
\EqTag{NKstruct}
\end{equation}
The functions $\mu _1,\mu_2  $ are the Monge-Amp\`ere invariants. 
\end{Theorem}
The singular or characteristic Pfaffian system for $\CalI$ are given by
\begin{equation}
\hV = \{\theta^i, \omega^1,\pi_1\}_{0\leq i \leq 2}\ , \quad \cV = \{ \theta^i,\omega^2,\pi_2\}_{0\leq i \leq 2} .
\EqTag{MAS2a}
\end{equation}
For the next lemma we note that the Cauchy characteristics for $\CalI'=\langle \theta^0\rangle_{\text diff}$ are given in the dual frame by
\begin{equation}
\Cau(\CalI') = \spn \{ \partial_{\pi^1}, \partial_{\pi^2} \} .
\EqTag{CIp}
\end{equation}
We assume $\Cau(\CalI')$ to be regular and $\bfq_{\Cau_{\CalI'}}:M \to M/\Cau_{\CalI'}$ to be the smooth quotient map.

The vanishing of the Monge-Amp\`ere invariants permits the following reduction lemma.

\begin{Lemma} \StTag{MRed} Suppose $\CalI$ is the standard rank 3 Pfaffian system for a hyperbolic PDE in the plane, and that the invariant conditions $\mu_1=\mu_2=0$ in equation \EqRef{NKstruct} are satisfied. Then
$\bar \CalI = \CalI/\bfq_{\Cau_{\CalI'}}$ is a hyperbolic Monge-Amp\`ere system on the five dimensional manifold $ M/\Cau_{\CalI'}$,  and $\bar\CalI =\CalI/\bfq_{\Cau_{\CalI'}}$ is the de-prolongation of $\CalI$.
\end{Lemma}
\begin{proof} We follow Section 4.2. The structure equations \EqRef{NKstruct} with $\mu_1=\mu_2=0$ determine the $\bfq_{\Cau_{\CalI'}}$ semi-basic forms as
\begin{equation}
\langle \CalI_{\bfq_{\Cau_{\CalI'}},sb}\rangle = \langle \theta^0, \omega^1\wedge \theta^1, \omega^2 \wedge \theta^2  \rangle_{\text alg} .
\EqTag{RIM}
\end{equation}
The quotient 
$\CalI/\bfq_{\Cau_{\CalI'}}=\CalI_{\bfq_{\Cau_{\CalI'}},sb} /\bfq_{\Cau_{\CalI'}}$, and  by equation \EqRef{RIM}  is a hyperbolic Monge-Amp\`ere EDS.   Finally prolongation using $\omega^1 \wedge \omega^2$ as independence condition proves the final claim. 
\end{proof}

\begin{Remark} The PDE $F=0$ satisfies the invariant conditions $\mu_1=\mu_2=0$, if and only if $F=0$ is a Monge-Amp\`ere equation, see \cite{gardner-kamran:1993a}:
\end{Remark}

The EDS $\CalI$ is said to be Darboux integrable (and not Monge integrable) if $\hV^\infty$ and $\cV^\infty$ are rank two, in which case we may assume equations \EqRef{NKstruct} and \EqRef{MAS2a} hold and
\begin{equation}
\hV^\infty=\{ \omega^1, \pi_1 \}\, \quad \cV^\infty = \{ \omega^2, \pi_ 2 \} .
\EqTag{MAS3a}
\end{equation}
We now recall a pivotal result  from   \cite{anderson-fels-vassiliou:2009a} on the geometric properties of Darboux integrable systems applied to the case of hyperbolic PDE in the plane.

\begin{Theorem} \StTag{Qrep} (The quotient representation of a Darboux integrable system)
 Let $\CalI$ be Darboux integrable with singular Pfaffian systems \EqRef{MAS2a} satisfying \EqRef{MAS3a}, and  let $x\in M$. Let $M_1$ be the five dimensional maximal integral manifold of $\hV^\infty$ through $x$ and let $\CalK_1$ be rank 3 Pfaffian system which is the restriction of $\hCalV$ to $M_1$. Similarly let $M_2$ be the five dimensional maximal integral manifold of $\cV^\infty$ though $x$ and $\CalK_2$ the rank 3 Pfaffian system which is the restriction of $\cCalV$ to $M_2$.  Then 
 
\noindent
{\bf [i]} there exists a pair of isomorphic three dimensional Lie algebra of vector-fields $\Gamma_a$ on $M_a$ each of which is point-wise linearly independent, transverse to $\CalK_a$ and symmetries of $\CalK_a$;

\noindent
{\bf [ii]}  there exists  open sets $U\subset  M$, $U_1\in M_1$ and $U_2\in M_2$ containing $x$ such that
\begin{equation}
\CalI|_U = (\CalK_1+\CalK_2)|_{U_1\times U_2}/\Gamma_{\text diag} .
\EqTag{EQrep}
\end{equation}
\end{Theorem}

The algebra $\Gamma_a$ is called the Vessiot algebra of $\CalI$ which is determined by a sequence of coframe adaptations, is a fundamental
invariant of Darboux integrable systems whose construction is given in \cite{anderson-fels-vassiliou:2009a} (see also \cite{anderson-fels:2011a}).

Before proceeding to the proofs of Theorems \StRef{DI5}, \StRef{BTMA0}, and \StRef{BTMA} we need to examine the derived sequence of $\CalK_a$. Using the inclusion map $i_1:M_1 \to U$
we have   $K^1_1 = i_1^*( \hV)$ which from \EqRef{MAS2a} is
to $M_1$ 
$$
K^1_1 =\spn \{ \ \theta^0,\ \theta^1,\ \theta^2\ \}\ .
$$
The structure equations for the Pfaffian system $\CalK_1$ are easily obtained by restricting equations \EqRef{NKstruct} to $M_1$. On $M_1$ we have $\omega^1=0$ and $\pi_1=0$ on $M_1$ and the structure equations from \EqRef{NKstruct} pullback to (by an abuse of notation),
\begin{equation}
\begin{aligned}
d\theta^0 & = \theta^2\wedge \omega^2  && \mod \theta^0 \\
d \theta^1 & = \mu_1 \theta^2 \wedge \pi^2 & & \mod \theta^0,\theta^1 \\
d \theta^2 & = -\pi_2 \wedge \omega^2 & & \mod \ \theta^0,\theta^1 .
\end{aligned}
\EqTag{MAS0}
\end{equation}
Analogous equations hold for $\CalK_2$ on $M_2$. These equations easily imply the following.
\begin{Corollary}\StTag{DSK} If the Monge-Amp\`ere invariants $\mu_1,\mu_2$ are nowhere vanishing then by equation \EqRef{MAS0} the derived series for $K_1$ and $K_2$ are
$$
(K_1)'=\{ \ \theta^0,\ \theta^1 \} ,\ (K_1)''=\{\, 0 \,\}; \quad (K_2)'=\{ \ \theta^0, \ \theta^ 2 \} , \ (K_2)''=\{\, 0 \, \} .
$$
If the Monge-Amp\`ere invariants satisfy $\mu_1=\mu_2=0$ then
$$
(K_1)'=\{ \ \theta^0,\ \theta^1 \} , \  (K_1)''=\{ \, \theta^1  \}, \ (K_1)'''=\{\, 0 \, \} ;
\quad (K_2)'=\{ \, \theta^0, \, \theta^ 2\, \} ,  (K_2)'=\{ \, \theta^2 \, \}, \ (K_2)'''=\{\, 0 \,  \} .
$$
\end{Corollary}

\subsection{B\"acklund Transformations for PDE in the plane: Proof of Theorem \StRef{DI5}}


\begin{proof} (Proof of Theorem \StRef{DI5})  Let $\CalI_2$ be the standard rank three Pfaffian system on a seven manifold $M$ representing a Darboux integrable hyperbolic PDE in the plane satisfying \EqRef{MAS3a}.  

We then construct diagram \EqRef{Intro10} as a commutative diagram of EDS by identifying  the EDS $\CalI$ in Theorem \StRef{ThIntro1} along with the three Lie algebras $\Gamma_H$, $\Gamma_{G_1}$ and $\Gamma_{G_2}$ with $\Gamma_H = \Gamma_{G_1} \cap \Gamma_{G_2}$. Applying Theorem \StRef{ThIntro1} produces diagram \EqRef{Intro10}.

Our starting point is Theorem \StRef{Qrep}. Let $x\in M$, $U,U_1,U_2$, $\Gamma_1$ and $\Gamma_2$ be 
as in Theorem \StRef{Qrep}. Since $\Gamma_1$ and $\Gamma_2$ are isomorphic, choose basis so that  $\Gamma_1={\text span}\{X_i\}_{1\leq i \leq 3}$, and $\Gamma_2={\text span} \{Y_i\}_{1\leq i \leq 3}$ have the same structure constants.  The infinitesimal action of $\Gamma_{G_2}=\Gamma_{\text diag}$ is
\begin{equation}
\Gamma_{G_2} = \Gamma_{\text diag}=\spn\{ X_i + Y_i\}_{1\leq i \leq 3},
\EqTag{GamG2}
\end{equation}
and equation \EqRef{EQrep} produces the right hand side of \EqRef{Intro10} where $\CalI= \CalK_1+\CalK_2$, is on $M_1\times M_2$ and $\Gamma_{G_2}=\Gamma_{\text diag}:U_1\times U_2 \to U$ as described in Theorem \StRef{Qrep}. We now identify $\Gamma_H$ and $\Gamma_{G_1}$.

A cursory  glance at the classification of
real 3-dimensional Lie algebras shows that  except for the case $\Gamma_1=\Gamma_2 = {\bf so}(3)$,  every  such algebra admits a 2-dimensional subalgebra $\lieh$. Therefore a basis for the infinitesimal generators $\Gamma_a$ may be chosen so that $X_1,X_2$ and $Y_1,Y_2$ are the infinitesimal generators for the action of $\lieh$ so that the structure equations satisfy,
\begin{equation}
	[X_1, X_2 ] =  \epsilon X_1 \quad\text{and}\quad [Y_1, Y_2 ] =  \epsilon Y_1 , \quad \epsilon = 0, 1.
\EqTag{Hadapt}
\end{equation}
In a neighbourhood of $x$, the sub-algebras $\{ X_1, X_2 \}$ and $\{ Y_1 , Y_2\}$ maybe chosen transverse to the two-dimensional derived systems $\CalK_1'$ and $\CalK_2'$. Relabel the corresponding open sets, $U_1,U_2,U$.


Using the notation in equation \EqRef{Hadapt}, let  $\Gamma_H \subset \Gamma_{G_2}$ be the two dimensional subalgeba  $\Gamma_{H} = \{ X_1+Y_1, X_2+Y_2\} $. Let $\bfq_{\Gamma_H}:U_1\times U_2 \to (U_1\times U_2)/\Gamma_H$ be the quotient map, and $\bfp_2$ the induced orbit map from Theorem  \StRef{ORD} with the subalgebra $\Gamma_H \subset \Gamma_{\text diag}$. Letting $\CalB = (\CalK_1+\CalK_2)/\Gamma_H$, we have completed the right triangle in diagram in \EqRef{Intro10}. Theorem \StRef{ORD} also shows that each of the elements in the right triangle of \EqRef{Intro10} is an integrable extension.

We now identify $\Gamma_{G_1}$. Using the basis of infinitesimal generators $\Gamma_a$ adapted to the two dimensional sub-algebras in \EqRef{Hadapt} let
\begin{equation}
\Gamma_{G_1} = \{ X_1, Y_1, X_2+ Y_2\}.
\EqTag{GamG1}
\end{equation}
Note that $\Gamma_H= \Gamma_{G_2} \cap \Gamma_{ G_1}$. Let $\bfq_{\Gamma_{G_1}}:U_1\times U_2 \to (U_1\times U_2)/\Gamma_{G_1}$ and let $\bfp_1$ the induced orbit map from Theorem \StRef{ORD} with the subalgebra $\Gamma_H \subset \Gamma_{G_1}$. Let $\CalI_{1} = (\CalK_1+\CalK_2)/\Gamma_{G_1}$. We can now use part {\bf [ii]} of Theorem \StRef{ThIntro1} to produced diagram \EqRef{Intro10} as
a commutative diagram of integral extensions. The fact that $\Gamma_{G_1} $ is transverse to $\CalK_1+\CalK_2$ follows from the fact
that $\Gamma_a$ are transverse to $\CalK_a$ (see Theorem \StRef{Qrep}).  This shows $\CalB=(\CalK_1+\CalK_2)/\Gamma_H$ is a B\"acklund transformation with one-dimensional fibres.

For the final part of the proof we show  $\CalI_1=(\CalK_1+\CalK_2)/\Gamma_{G_1}$ is locally equivalent to the EDS $\CalI_{s=0}$ which is the standard rank three Pfaffian system  on a seven manifold for the wave equation. To do this we use Theorem  6.1 in \cite{anderson-fels:2011a} (see also Lemma 6.7 in \cite{anderson-fels:2011a}). Let $\tilde \Gamma_1=\{X_1\}, \tilde \Gamma_2=\{Y_1\}$, then the reference above,
\begin{equation}
\CalI_1=(\CalK_1+\CalK_2)/\Gamma_{G_1}= (\CalK_1/\tilde\Gamma_1 + \CalK_2/\tilde \Gamma_2)/ (\Gamma_{G_1}/(\tilde \Gamma_1 \times \tilde \Gamma_2))
\EqTag{DQ}
\end{equation}
Now by the choice of $\{X_1,X_2\}$ and $\{Y_1, Y_2\}$ being transverse to $\CalK_a'$, the quotient system
$\CalK_1/\tilde\Gamma_1 + \CalK_2/\tilde \Gamma_{2} $ is a Pfaffian system (see \cite{anderson-fels:2005a}) where each $\CalK_a/\tilde \Gamma_a$ has derived sequence rank $[2,1,0]$. 
 Furthermore the quotient algebras $\{X_1,X_2\}/\{ X_1\} $ and $\{Y_1,Y_2\}/\{Y_1\}$ are transverse to $(\CalK_a')/\tilde \Gamma_a$.
The action of one-dimensional Lie algebra $\Gamma_{G_1}/ (\tilde \Gamma_1 \times \tilde \Gamma_2)=\{ X_2+Y_2\} $ is therefore transverse to $(\CalK_1'/\tilde \Gamma_1+\CalK_2'/\tilde \Gamma_2)$. Again using \cite{anderson-fels:2011a} we conclude from these remarks and equation \EqRef{DQ} that $(\CalK_1+\CalK_2)/\Gamma_{G_1}$ is a rank 3 hyperbolic Pfaffian system.

We use Theorem  6.1 part {\bf [i]} in \cite{anderson-fels:2011a} to compute number of Darboux invariants, which states $\CalI/\Gamma_{G_1}$ has $M_1/\pi_1(\Gamma_{G_1})$ 
and $M_2/\pi_2(\Gamma_{G_1})$ invariants where $\pi_a:M_1\times M_2 \to M_a $, a=1,2.
 We have from equation \EqRef{GamG1},
$$
\pi_1(\Gamma_{G_1}) = \spn\{ X_1, X_2 \}, \quad \pi_2(\Gamma_{G_1}) = \spn \{ Y_1, Y_2\}.
$$
Since $M_1$ and $M_2$ are 5 dimensional, we conclude (equation 6.10 in \cite{anderson-fels:2011a}) that
$$
\rank \hV_0^\infty =\dim M_2-\rank \pi_2(\Gamma_{G_1})=  3, \quad \rank \cV_0^\infty = \dim M_1 - \rank \pi_1(\Gamma_{G_1}) = 3.
$$

If we now show that $\CalI_1$ is the EDS for a PDE in the plane, the local equivalence with $\CalI_{s=0}$ will follow from a well-known theorem of Lie's which states: a hyperbolic PDE in the plane where each characteristic system admits three Darboux invariants is locally equivalent to $\CalI_{s=0}$. 

By the transversality of the one-dimensional Lie algebra $\Gamma_{G_1}/ (\tilde \Gamma_1 \times \tilde \Gamma_2)=\{ X_2+Y_2\} $ to the rank two distribution $\CalK_1'/\tilde \Gamma_1+\CalK_2'/\tilde \Gamma_2$
we have $\rank \CalI_1'=1$. Let $\theta $ be a local generator for $\CalI'$. Then $\CalI_1$ represents a PDE in the plane if and only if $\theta$ has rank 2 (see \cite{stormark:2000a}).  Let $\theta_{sb} \in \CalK_1'/\tilde \Gamma_1+\CalK_2'/\tilde \Gamma_2$, be $\Gamma_{G_1}/ (\tilde \Gamma_1 \times \tilde \Gamma_2)=\{ X_2+Y_2\} $ semi-basic. Then $\theta$ has rank 2 if and only if $\theta_{sb}$ is rank 2 because $\bfq^*_{\Gamma_{G_1}} \theta = \lambda \theta_{sb}$.  

By the transversality of the action of $\{X_1,X_2\}/\tilde \Gamma_1$ to $\CalK_1'/\tilde \Gamma_1$ on $M_1/\tilde \Gamma_1$ and $\{ Y_1, Y_2\}/\tilde \Gamma_2$
to  $\CalK_2'/\tilde \Gamma_2$ on $M_2/\tilde \Gamma_2$ we may choose $\theta^1 \in \CalK_1'/\tilde \Gamma_1$ and $\eta^1 \in \CalK_2'/\tilde \Gamma_2$ satisfying $\theta^1(X_2) = 1, \eta^1(Y_2)= 1$. Then
$$
\theta_{sb} = \theta^1 - \eta^1 
$$
is semi-basic. The derived system for $\CalK_a/\tilde \Gamma_a$ are $[2,1,0]$ so we may choose one-forms $\omega^1,\theta^2$ on $M_1/\tilde \Gamma_1$ and $\sigma^1,\eta^2$ on $M_2/\tilde \Gamma_2$ so that $\CalK_1/\tilde \Gamma_1 =\langle \theta^1, \theta^2\rangle_{\text diff}$, $\CalK_2/\tilde \Gamma_2 =\langle \eta^1, \eta^2\rangle_{\text diff}$ and
$$
d \theta^1= \sigma^1 \wedge \theta^2\quad \mod \theta^1 ,\quad d \eta^1 = \tau^1 \wedge \eta^2 \quad \mod  \eta^1.
$$
Therefore $\theta^1_{sb} \wedge (d \theta_{sb}^1)^2 \neq 0$ and so both $\theta_{sb}$ and  $\theta$ have rank 2. \end{proof}

\begin{Example} The non-linear partial differential equation $3u_{xx}(u_{yy})^3 +1 =0$ is Darboux integrable. Applying the algorithm in \cite{anderson-fels:2011a} produces the following quotient representation. Let $\CalK_1+\CalK_2$ be the sum of two copies of the standard Pfaffian system for the Cartan-Hilbert equation $u'=(v'')^2$, on $\real^5\times \real^5$ with coordinates $(s,u,v,v_s,v_{ss}; t, y,z, z_t,z_{tt})$,
\begin{equation}
\CalK_1+\CalK_2 =\langle du -(v_{ss})^2 ds , dv- v_s ds, dv_s-v_{ss}ds, dy -(z_{tt})^2dt,dz-z_t dt, dz_t-z_{tt}dt \rangle_{\text diff}\, .
\end{equation}
The Lie algebra $\Gamma_{G_2} $ in Theorem \StRef{Qrep} is
$$
\Gamma_{G_2} = \spn \{ \partial_v-\partial_z,  s\partial_v + \partial_{v'}+ t\partial_z + \partial_{z'} , \partial_u-\partial_y \} .
$$
Now let
\begin{equation}
\begin{aligned}
\Gamma_{G_1}  &= \spn \{ \partial_v, \partial_z,  s\partial_v + \partial_{v_s}+ t\partial_z + \partial_{z_t}  \} \\
\Gamma_{H}  &=  \Gamma_{G_1} \cap \Gamma_{G_2} = \spn \{ \partial_v-\partial_z,  s\partial_v + \partial_{v_s}+ t\partial_z + \partial_{z_t}   \}  .
\end{aligned}
\end{equation}
Coordinates can be chosen on the 8-manifold $M/\Gamma_H$ so that
the quotient map $\bfq_H:M \to M/\Gamma_H$ is
\begin{equation}
\bfq_H =        \left( x_1 = t,\  x_2 = s,\  x_3 = z_{tt},\  x_4 = v_{ss},\  x_5 = z_s -  v_s,\  x_6 = u,\  x_7 = y, \ x_8 = z + v -(s+t)v_s \right)\, .
\end{equation}
Coordinates on $M/\Gamma_{G_1}$ can be chose so the the map $\bfp_1$ is
\begin{equation}
\begin{gathered}
\bfp_1=	\left(\ X' =- x_6  + x_2 x_4^2,   \quad    Y' = x_7 - x_1 x_3^3,\quad  U'=   x_5  - x_1x_3 + x_2x_4,\   \right.
\\ 	
	\left. \ P'= \frac1{2x_4} ,\quad Q' = \frac1{2x_3}, \quad      R' = -\frac{1}{2x_2(x_4)^2}, \quad  T' = \frac{1}{2x_1(x_3)^2} \right) .
\end{gathered}
\EqTag{ExB8}
\end{equation}
Finally coordinates on $M/\Gamma_{G_2}$ can be chosen so that $\bfp_2$ has the form
\begin{equation}
\begin{aligned}
\bfp_2=\left(
X  = -2 \frac{x_3 +x_4}{x_1+x_2},  \quad Y =x_5 - \frac12 (x_3-x_4)(x_1+x_2),\quad   S =   \frac12(x_2- x_1), \quad T= \frac2{x_1+x_2},\right.
\\
	U  =  2 \frac{x_3 +x_4}{x_1 +x_2} (-x_8  +x_1x_5) - x_6   +\frac13(2x_1 -x_2)x_4^2
	-\frac23(x_1 +x_2) x_3x_4 -\frac13(x_1 - 2x_2)x_3^2,
	\\
\left.	P  =  x_8 - x_1x_5 + \frac16 (x_1 +x_2)((2x_1 -x_2)x_3  -  (x_1-2x_2)x_4) ,\quad 
	Q  = \frac12\frac{x_4x_1-x_3x_2}{x_1+x_2} \right) .
\end{aligned}
\EqTag{ExB10}
\end{equation}
The choice of variables \EqRef{ExB8} and \EqRef{ExB10}  give the B\"acklund transformation in diagram \EqRef{Intro10}
 as
\begin{equation}
\begin{gathered}
\beginDC{\commdiag}[3]
\Obj (1,28)[B]{$\boxed{\begin{gathered}  dx_5 - x_3\,dx_1 + x_4\,dx_2 \\[-2\jot]  dx_6- x_4^2\,dx_2 ,\    dx_7 - x_3^2\, dx_1  \\[-2\jot]    dx_8 - x_5\,dx_1 +x_4(x_1+x_2)\,dx_2     \end{gathered}       } $}
\Obj(-28,0)[I1]{$\boxed{\begin{gathered}  dU'- P' dX'-Q'dY' \\[-2\jot]  dP'-R'dX',\    dQ'-T'dY'   \end{gathered}  }$}
\Obj(35, 0)[I2]{$\boxed{   \begin{gathered}  dU- P dX-QdY \\[-2\jot]  dP+\frac{1}{3T^3}dX-SdY,\    dQ-SdX-TdY   \end{gathered}    }$ . }
\mor(2, 17)(-32,7){$\bfp_1$}[\atright, \solidarrow]
\mor(2,17)(38,7){$\bfp_2$}[\atleft, \solidarrow]
\enddc
\end{gathered}
\end{equation}
Where the standard representation the for PDE $U'_{X'Y'} = 0$ and $3U_{XX}U_{YY}^3 +1 =  0$ as EDS appear.

\end{Example}

\subsection{The Monge-Amp\`ere Case: Proof of Theorems \StRef{BTMA0} and \StRef{BTMA}}

Let $\CalI_2 $ be a rank 3 Pfaffian system on a 7 manifold $M$ for a Darboux integrable hyperbolic Monge-Amp\`ere equation satisfying \EqRef{MAS3a}. Let $ \bar \CalI_2=\CalI_2/\bfq_{\Cau_{\CalI_2'}}$ be the corresponding Monge-Amp\`ere system on the five-manifold $M/\Cau_{\CalI_2'}$ (Lemma \StRef{MRed}), where $\bfq_{\Cau_{\CalI_2'}}:M \to M/\Cau_{\CalI_2'}$ is the smooth quotient map. Likewise let $\Cau(\CalI'_{s=0})$ be
the two dimensional space of Cauchy characteristics for  the derived system $\CalI'_{s=0}$ with smooth quotient map 
$\bfq_{\Cau_{\CalI_{s=0}'}}:M_{s=0} \to M_{s=0}/\Cau_{\CalI_{s=0}'}$. The EDS $\bar \CalI_{s=0}= \CalI_{s=0}/\bfq_{\Cau_{\CalI'_{s=0}}}$ is the Monge-Amp\`ere form of the wave equation.

\begin{proof}  (Theorem \StRef{BTMA0}) Let $x\in M$ and  recall that  $M_1$ is the maximal integral manifold of $\hV^\infty$ through $x$  and $\CalK_1$ is the restriction of $\hCalV$ to $M_1$ with the analogous definition of $M_2$ and $\CalK_2$ (Theorem \StRef{Qrep}).
By the proof of Theorem \StRef{DI5} in Section 5.2 there exists open sets $U\subset M,U_1\subset M_1,U_2\subset M_2$ where  $x\in U$ and diagram \EqRef{Intro10} holds. In order to produce diagram \EqRef{FCC} we need to de-prolong diagram \EqRef{Intro10} by applying Theorem \StRef{SymmodC}.

By Corollary \StRef{DSK} we have that $\CalK_1'=\{ \theta^0, \theta^1\}$, and from equation \EqRef{MAS0} with $\mu_1=0$,  that $\spn\{ \partial_{\pi_2}\} $ is the Cauchy characteristic of $\hCalK_1'$. A similar argument holds for $\CalK_2$ on $U_2$ and so the Cauchy characteristics for $\CalK_1'+\CalK_2'$ on $U_1\times U_2$ are
$$
\Cau(\CalK_1'+\CalK_2') |_{U_1\times U_2}=\spn \{ \partial_{\pi_2}+0, 0+\partial_{\pi_1} \} .
$$
By part {\bf [ii]} in Theorem \StRef{SymmodC} let  $\tilde \Gamma_{G_a}$ and $\tilde \Gamma_H$
be the projected Lie algebras of symmetries of $(\CalK_1+\CalK_2)/\bfq_{\Cau_{\CalK_1'+\CalK_2'}}$ from $\Gamma_{G_a}$ and $\Gamma_H$ on $M_1\times M_2$. 
In the proof of Theorem \StRef{DI5} in Section 5.2 it was noted that all the algebras $\Gamma_{G_a}$ and $\Gamma_H$ 
are transverse to $\CalK_1'+\CalK_2'$.  This combined with the fact  $\Cau(\CalK_1'+\CalK_2') \subset \ann(K_1^1+K_2^1)$ implies  $\Cau(\CalK_1'+\CalK_2')/\Gamma_H$ and $\Cau(\CalK_1'+\CalK_2')/\Gamma_{G_a}$ are rank 2. Therefore part {\bf [v]} of  Theorem \StRef{SymmodC} and dimensional considerations 
imply 
\begin{equation}
\Cau(\CalB') =\Cau(\CalK_1'+\CalK_2')/\Gamma_H, \quad \Cau(\CalI'_2) = \Cau(\CalK_1'+\CalK_2')/\Gamma_{G_2}, \quad \Cau(\CalI_{s=0}')=\Cau(\CalK_1'+\CalK_2')/\Gamma_{G_1}.
\EqTag{CAUR}
\end{equation}
Diagram \EqRef{Ddepro} (or Remark \StRef{CPD}) along with equation \EqRef{CAUR} then give the de-prolongation diagrams,
\begin{equation}
\begin{gathered}
\beginDC{\commdiag}[3]
\Obj(0, 20)[I1]{$\CalK_1+\CalK_2$}
\Obj(0, 0)[I1G]{$(\CalK_1+\CalK_2)/\bfq_{\Cau_{\CalK_1'+\CalK_2'}} $}
\Obj(35, 20)[I1C]{$\CalB $}
\Obj(35, 0)[I1Z]{$\bar \CalB$}
\mor{I1}{I1G}{$\bfq_{\Cau_{\Cau_{\CalK_1'+\CalK_2'}} }$}[\atleft, \solidarrow]
\mor{I1G}{I1Z}{$\tilde \bfq_{\Gamma_H}$}[\atleft, \solidarrow]
\mor{I1C}{I1Z}{$\bfq_{\Cau_{\CalB'}}$}[\atright, \solidarrow]
\mor{I1}{I1C}{$\bfq_{\Gamma_H}$}[\atright, \solidarrow]
\enddc
\end{gathered}
\qquad
\begin{gathered}
\beginDC{\commdiag}[3]
\Obj(0, 20)[I1]{$\CalK_1+\CalK_2$}
\Obj(0, 0)[I1G]{$(\CalK_1+\CalK_2)/\bfq_{\Cau_{\CalK_1'+\CalK_2'} }$}
\Obj(35, 20)[I1C]{$\CalI_2 $}
\Obj(35, 0)[I1Z]{$\bar \CalI_2$}
\mor{I1}{I1G}{$\bfq_{\Cau_{\Cau_{\CalK_1'+\CalK_2'}} }$}[\atleft, \solidarrow]
\mor{I1G}{I1Z}{$\tilde \bfq_{\Gamma_{G_2}}$}[\atleft, \solidarrow]
\mor{I1C}{I1Z}{$\bfq_{\Cau_{\CalI_2'}}$}[\atright, \solidarrow]
\mor{I1}{I1C}{$\bfq_{\Gamma_{G_2}}$}[\atright, \solidarrow]
\enddc
\end{gathered}
\EqTag{Ddepro12}
\end{equation}
with a similar diagram holding for $\Gamma_{G_1}$ and $\CalI_{s=0}$.  Now applying Theorem  \StRef{ThIntro1} with $\tilde \Gamma_H =\tilde \Gamma_{G_1} \cap \tilde \Gamma_{G_2}$  produces the diagram
\begin{equation}
\beginDC{\commdiag}[3]
\Obj(0, 30)[I]{$(\CalK_1+\CalK_2)/\bfq_{\Cau_{\CalK_1'+\CalK_2'}} $}
\Obj(0, 10)[H]{$\bar \CalB$}
\Obj(-30, 0)[I1]{$\bar \CalI_{s=0}$}
\Obj(30,  0)[I2]{$\bar \CalI_2$}
\mor{I}{H}{$\tilde \bfq_{\Gamma_H}$}[\atleft, \solidarrow]
\mor{H}{I1}{$\tilde \bfp_1$}[\atleft, \solidarrow]
\mor{H}{I2}{$\tilde \bfp_2$}[\atright, \solidarrow]
\mor{I}{I1}{$\tilde \bfq_{\Gamma_{G_1}}$}[\atright, \solidarrow]
\mor{I}{I2}{$\tilde \bfq_{\Gamma_{G_2}}$}[\atleft, \solidarrow]
\enddc
.
\EqTag{Ddepro13}
\end{equation}
Combining the bottom edges of diagrams \EqRef{Ddepro12} together with \EqRef{Ddepro13} produces the bottom part of diagram \EqRef{FCC}.
\end{proof}

\begin{proof} (Theorem \StRef{BTMA}).  Theorem \StRef{BTMA} will follow immediately from diagram \EqRef{Ddepro13} in  the previous proof once we identify $\CalK_1$ and $\CalK_2$. By Corollary \StRef{DSK}, $\CalK_a$ have derived flag $[3,2,1,0]$. Therefore Engle's Theorem \cite{bryant-chern-gardner-griffiths-goldschmidt:1991a}  applies and $\CalK_a$  may be identified locally with the canonical contact system on $J^3(\real,\real)$. Finally $\CalC_a = \CalK_a/\bfq_{\Cau_{\CalK_a'}}$ are the de-prolongation of $\CalK_a$, and may be identified  locally with the canonical contact system on $J^2(\real,\real)$. This finishes the proof of Theorem \StRef{BTMA} and produces diagram \EqRef{Intro1b}
\end{proof}




\begin{bibdiv}
\begin{biblist}





\bib{anderson-fels:2005a}{article}{
  author={Anderson, I. M.},
  author={Fels, M. E.},
  title={Exterior Differential Systems with Symmetry},
  journal={Acta. Appl. Math.},
  year={2005},
  volume={87},
  pages={3--31},
}


\bib{anderson-fels:2011a}{article}{
  author={Anderson, I. M.},
  author={Fels, M. E.},
  title={B\"acklund Transformations for Darboux Integrable Differential Systems},
  journal={arXiv:1108.5443v1},
  year={2011}
}


\bib{anderson-fels-vassiliou:2009a}{article}{
  author = {Anderson, I. M.},
  author = {Fels, M. E.},
  author = {Vassiliou, P. V.}, 	
  title={Superposition Formulas for Exterior Differential Systems},
  journal ={Advances in Math.},
  volume = {221},
  year = {2009},
  pages = {1910 -1963}
}


\bib{bryant-chern-gardner-griffiths-goldschmidt:1991a}{article}{
  author = {Bryant, R. L.},
  author = {Chern, S. S.},
  author = {Gardner, R.B.},
  author = {Goldschmidt, H. L.},	
  author = {Griffiths, P. A.},
   title = {Exterior Bifferential Systems},
 series = { MSRI Publications},
  year={1991},
  volume={18}
}

\bib{bryant-griffiths:1995a}{article}{
  author={Bryant, R. L},
  author={Griffiths, P. A.},
  title={Characteristic cohomology of differential systems II: Conservation Laws for a class of Parabolic equations},
  journal={Duke Math. Journal},
  volume={78},
  number={3},
  year={1995},
  pages={531--676},
}


%





\bib{Clelland-Ivey:2009a}{article}{
 author = {Clelland, J. N.},
 author = {Ivey, T. A.},
 title = {B\"acklund transformations and Darboux integrability for nonlinear wave equations},
 volume = {13},
 journal = {Asian J. Math},
 pages={15--64}
}


\bib{fels:2007a}{article}{
  author = {Fels, M. E.}, 	
  title={Integrating ordinary differential equations with symmetry revisited},
  journal ={Foundations of Computational Math},
  volume = {7({\bf 4})},
  year = {2007},
  pages = {417-454}
}


\bib{gardner-kamran:1993a}{article}{
author = {R. B. Gardner},
author = { N. Kamran},
title = {Characteristics and the geometry of hyperbolic equations in the plane},
year =  {1993},
volume = {104},
journal = { J. Differential Equations},
pages = { 60--117}
}






\bib{ivey-langsberg:2003}{book}{
 author = {T. A. Ivey},
 author = {J. M. Landsberg},
 title = {Cartan for Beginners: Differential Geometry via Moving Frames and Exterior Differential Systems},
series = {Graduate Studies in Mathematics},
volume = {61},
publisher = {Amer. Math. Soc.}, year={2003}
}


\bib{munkres:2000a}{book}{
  author={Munkres, J.R.},
  title={Topology, Second Edition},
  year={2000},
  publisher={Prentice Hall},
  address={USA},
}





\bib{stormark:2000a}{article}{
  author={Stormark, O.},
  title={Lie's structural approach to PDE systems},
  series={Encyclopedia of Mathematics and its Applications},
  volume={80},
  publisher={Cambridge Univ. Press},
  address={Cambridge, UK},
  year={2000},
}



	

\end{biblist}
\end{bibdiv}
\end{document}